\documentclass[12pt]{amsart}

\usepackage{amssymb}
\usepackage{amsmath}
\usepackage[english]{babel}
\usepackage{amsfonts}

\newtheorem{thm}{Theorem}

\newtheorem{lemma}{Lemma}
\newtheorem{cor}{Corollary}
\newtheorem{prop}{Proposition}

\theoremstyle{definition}
\newtheorem{defn}{Definition}
\newtheorem{remk}{Remark}
\newtheorem{exam}{Example}

\newtheorem*{thm*}{Theorem}
\newtheorem*{prop*}{Proposition}

\newcommand{\cyc}[1]{\mathbb{Z}/{#1}}
\newcommand{\Un}[1]{\mathbf{1}_{#1}}
\newcommand{\set}[1]{\left\{#1\right\}}
\newcommand{\diag}[1]{\mathrm{diag}\left(#1\right)}

\newcommand{\cA}{\mathcal{A}}
\newcommand{\cB}{{\mathcal{B}}}

\newcommand{\cD}{{\mathcal{D}}}
\newcommand{\cE}{{\mathcal{E}}}

\newcommand{\cC}{{\mathcal{C}}}
\newcommand{\cF}{{\mathcal{F}}}
\newcommand{\cH}{{\mathcal{H}}}
\newcommand{\cL}{{\mathcal{L}}}

\newcommand{\cR}{{\mathcal{R}}}
\newcommand{\cS}{{\mathcal{S}}}

\newcommand{\cV}{{\mathcal{V}}}
\newcommand{\cW}{{\mathcal{W}}}
\newcommand{\cZ}{{\mathcal{Z}}}

\newcommand{\bR}{{\mathbf{R}}}
\newcommand{\bC}{{\mathbf{C}}}
\newcommand{\bK}{{\mathbf{K}}}
\newcommand{\bL}{{\mathbf{L}}}

\newcommand{\dC}{{\mathbb{C}}}

\newcommand{\dL}{{\mathbb{L}}}

\newcommand{\dN}{{\mathbb{N}}}
\newcommand{\dR}{{\mathbb{R}}}
\newcommand{\dT}{{\mathbb{T}}}
\newcommand{\dZ}{{\mathbb{Z}}}

\newcommand{\ga}{\mathfrak{a}}
\newcommand{\gb}{\mathfrak{b}}

\newcommand{\gee}{\mathfrak{e}}
\newcommand{\gf}{\mathfrak{f}}

\newcommand{\gy}{\mathfrak{y}}
\newcommand{\gx}{\mathfrak{x}}
\newcommand{\gz}{\mathfrak{z}}

\newcommand{\gog}{\mathfrak{g}}

\newcommand{\gp}{\mathfrak{p}}
\newcommand{\gP}{\mathfrak{P}}

\newcommand{\gA}{\mathfrak{A}}

\newcommand{\ov}{\overline}

\newcommand{\Aut}{{\mathrm{Aut}}}

\newcommand{\Id}{\mathrm{Id}}

\newcommand{\Rep}{\mathrm{Rep}}

\newcommand{\Quot}{\mathrm{Quot}}

\newcommand{\tr}{\mathrm{tr}}
\newcommand{\id}{\mathrm{id}}

\newcommand{\img}{\mathbf{i}}
\newcommand{\eps}{\varepsilon}

\title{Positivstellens\"atze for Algebras of Matrices}

\begin{document}
\author{Yurii Savchuk}
\address{Universit\"at Leipzig, Mathematisches Institut Johannisgasse 26, 04103 Leipzig, Germany}
\email{savchuk@math.uni-leipzig.de}
%\thanks{The first author was supported by the International Max Planck Research School for Mathematics in the Sciences (Leipzig)}

\author{Konrad Schm\"udgen}
\address{Universit\"at Leipzig, Mathematisches Institut Johannisgasse 26, 04103 Leipzig, Germany}
\email{schmuedgen@math.uni-leipzig.de}


\date{\today}


\maketitle
\begin{abstract}
The paper is concerned with various types of noncommutative Positivstellens\"atze for the matrix algebra $M_n(\cA)$, where $\cA$ is an algebra of operators acting on a unitary space, a path algebra, a cyclic algebra or a formally real field. Some new types of Positivstellens\"atze are proposed and proved. It is shown by examples that they occur. There are a number of results stating that a type of Positivstellensatz is valid for $M_n(\cA)$ provided that it holds for $\cA$.
\end{abstract}

\section{Introduction and preliminaries}
Positivstellens\"atze in real algebraic geometry express positive or non-negative polynomials on semi-algebraic sets in terms of weighted sums of squares of polynomials \cite{pd},\cite{mar}. They can be considered as generalizations of E. Artin's theorem on the solution of Hilbert 17th problem. Non-commutative Positivstellens\"atze are fundamental results of a new emerging mathematical field that might be called non-commutative real algebraic geometry (see \cite{s2} for a recent survey and some basic concepts). In the last decade a number of non-commutative Positivstellens\"atze have been obtained for various classes of $*$-algebras (see e.g. \cite{he}, \cite{hpm},\cite{sweyl}, \cite{senv}, \cite{sfrac}, \cite{cim2}). Despite of all these results it is not even clear what a proper generalization of Artin's theorem for non-commutative $*$-algebras should be (some proposals have been made in \cite{s2}, Section 4.1).

The purpose of the present paper is threefold. Our main aim is to prove non-commutative Positivstellens\"atze for $*$-algebras of matrices over various classes of (commutative or non-commutative) unital $*$-algebras or $*$-fields. The corresponding Positivstellens\"atze will be precisely the theorems stated in the sections of the paper. Secondly, our emphasize will be on different versions of such Positivstellens\"atze especially concerning the involved numerator and denominators sets. In the course of this we discover a number of new types of non-commutative Positivstellens\"atze. All of them can be viewed as generalizations of Artin's theorem to non-commutative $*$-algebras. And finally, we want to elaborate some methods and notions that can be used to develop non-commutative Positivstellens\"atze. Apart from diagonalization techniques for matrices the notion of a conditional expectation will play a crucial role in this respect.

In the various sections of this paper we are concerning with different $*$-algebras of matrices and different non-commutative Positivstellens\"atze. Let us explain, slightly symplifying, what all our Positivstellens\"atze have in common. There are always "natural" notions of positivity in the corresponding algebras. Often positive elements are those elements which act as positive operators on unitary spaces or have positive point evaluations. In many cases they are defined as elements which are positive with respect to all possible $*$-orderings. It can be shown that in many cases these positivity notions are in fact equivalent. We shall denote the positive elements of a $*$-algebra $\cA$ by $\cA_+$. A Positivstellensatz expresses, roughly speaking, positive elements in algebraic terms built on sums of squares by allowing denominator sets.

Let us briefly explain basic types of non-commutative Positivstellens\"atze that will appear in this paper. All of them can be considered as possible generalizations of Artin's theorem to general $*$-algebras. For let $\cA$ denote a unital $*$-algebras, $\cA^\circ$ the set of elemens sof $\cA$ which are not zero divisors and $\cA_+$ the set of positive elements of $\cA.$ We shall say that a
\begin{itemize}
\item \textit{Positivstellensatz of type I} holds for $x\in \cA_+$ if $x\in\sum\cA^2,$ 
\item \textit{Positivstellensatz of type II} holds for $x\in \cA_+$ if there exists a $c\in\cA^\circ$ such that $$c^*xc\in\sum\cA^2,$$
\item \textit{Positivstellensatz of type III} holds for $x\in \cA_+$ if exists an $c\in\cA^\circ\cap\sum\cA^2$ such that $$xc=cx \ \mbox{and}\ \ xc \in\sum\cA^2,$$
\item \textit{Positivstellensatz of type IV} holds for $x\in \cA_+$ if are pairwise commuting elements $c_1,\dots,c_k \in\cA^\circ\cap\sum\cA^2$ and pairwise commuting elements $b_1,\dots,b_n \in\sum\cA^2$ such that $$ xc_j=c_jx\ \mbox{for}\ \ j=1,\dots k\ \mbox{and}\ \ x c_1\dots c_k=b_1 \dots b_n.$$
\end{itemize}
If a Positivstellensatz is valid for all $x \in \cA_+$, we will simply say that it holds for $\cA$.

Probably the most general version of a Positivstellensatz is obtained if one allows general denominator sets and numerators sets as defined in Section \ref{sect_prelim}. There are further natural versions by adding conditions on the sets of denominators (for instance, by requiring denominator from the center or from distinguished commutative subalgebas). Also it should be emphasized that all these versions and types essentially depend on the chosen set of positive elements $\cA_+.$

Of course, type I is the most desirable version. There are quite a few $*$-algebras for which a Positivstellensatz of type I holds. Type I is valid for the rational functions $\dR(x_1,\dots,x_d)$ (by Artin's theorem), for the free polynomial algebra in $d$ generators (by Helton's theorem), by the trigonometric polynomials in one variable (by the Riesz-Fejer theorem) and for $*$-algebra generated by the shift operator (as shown recently in \cite{ss2}). Type I does not hold for the polynomial algebra $\dR[x_1,\dots,x_d]$, but type II does (again by Artin's theorem). The strict Positivstellens\"atze proved in \cite{sweyl,senv,sfrac} are all of type II. We shall see in Section \ref{subsect_ctrexam} below that there exist positive elements of a cyclic $*$-algebra for which a Positivstellensatz of type II is not valid, but there is a corresponding result with products of commuting squares in the denominator or likewise in the numerator. In the above terminology this means that a Positivstellens\"atze of types III and IV hold. This observation was in fact the starting point for our search to more general versions of Positivstellens\"atze. 

The following simple example illustrates how the product of a positive element with commutings sums of squares becomes a sum of squares.
\begin{exam}\label{exam_weyl_alg} Let $\cA$ be the Weyl algebra $\dC\langle
a,a^*|aa^*-a^*a=1\rangle$. The $*$-algebra $\cA$ acts on the unitary space $\cD$ of all finite complex sequences $\varphi=(\varphi_0,\dots, \varphi_n,0,\dots)$ with scalar product
$\langle \varphi,\psi\rangle= \varphi_0\ov{\psi}_0 + \varphi_1 \ov{\psi_1}+\varphi_2\ov{\psi}_2+\dots
$
%The $*$-algebra acts $\cA$ acts on $\cD$ 
by
$$
a\varphi=(\varphi_1,\sqrt{2}\varphi_2,\sqrt{3}\varphi_3,\dots)\ \mbox{and}\ \
a^*\varphi=(0,\varphi_1,\sqrt{2}\varphi_2,\sqrt{3}\varphi_3,\dots).
$$
Define $\cA_+=\{x \in \cA: \langle x\varphi,\varphi \rangle \geq 0 ~~{\rm for } ~~\varphi \in \cD\}$ and $N=a^*a.$ Then $N$ acts diagonal on the orthonormal basis $e_n:=(\delta_{nk})$, $n\in \dN_0$, of the unitary space $\cD$, that is, $Ne_n=ne_n$ for $n\in \dN_0$.
Hence a polynomial $f(N)\in \dC[N]$ is in $\cA_+$ if and only if $f(n)\geq 0$ for all $n\in \dN_0$. 

Since $aa^*-a^*a=1$, $c_k:=(a^*)^ka^k=N(N{-}1)\cdots (N{-}(k{-}1))\in \sum\cA^2.$ It can be proved (see \cite{fs} or \cite{ss}) that $f(N)\in\sum\cA^2$ if and only if there are polynomials $g_0,\dots,g_k\in\dC[N]$ such that
\begin{gather}\label{bsuma2}
f(N)=g_0(N)^*g_0(N)+c_1g_1(N)^*g_1(N)+\cdots+c_kg_k(N)^*g_k(N).
\end{gather}
%Since $aa^*-a^*a=1$, we have $c_k:=(a^*)^ka^k=N(N{-}1)\cdots (N{-}(k{-}1))\in \sum\cA^2$.
Set $x_n=(N{-}n)(N{-}(n{+}1))$ for $n\in \dN.$ Then we have $x_n\in \cA_+$ and $x_n\notin \sum \cA^2$, because $x_n$ is not of the form (\ref{bsuma2}). But $c_nx_n=x_nc_n$ and $c_nx_n =c_{n+2}\in \sum \cA^2$.

Let $0< n_1<n_2<\dots<n_m.$ Then $x:=x_{n_1}\cdots x_{n_m}\in \cA_+$ and $c_{n_1},\dots,c_{n_m}$ are pairwise commuting elements of $\sum \cA^2$ such that $xc_{n_j}=c_{n_j}x$ for $j{=}1,\dots,m$ and $c_{n_1} \cdots c_{n_m} x = c_{n_1+2}\cdots c_{n_m+2}.$ It is not diffcult to check that the set of polynomials of the form (\ref{bsuma2}) is closed under multiplication. Hence $c_{n_1}\cdots c_{n_m} x\in \cA^2$.
\end{exam}

Let us describe the content of this paper. Section \ref{sect_prelim} collects some basic definitions and facts which are used throughout the text. In a recent paper \cite{po} it was proved that a Positivstellensatz of type I is valid for path algebras. In Section \ref{sect_sos_path} we apply a conditional expectation from an appropriate matrix algebra to give a simple alternative proof of this result. Section \ref{sect_nc_matr} deals with the matrix $*$-algebra $M_n(\cA)$ over a unital $*$-algebra $\cA$ of operators acting on a unitary space such that $\cA$ has no zero divisors and $\cA^\circ:=\cA\backslash \{0\}$ is a left Ore set. By developing a non-commutative diagonalization procedure we prove that a Positivstellensatz of type II holds for $M_n(\cA)$ provided that it holds for $\cA$. Section \ref{sect_sos_cross} is concerning with the cross product algebra $\cA \times_\alpha G$ of a $*$-algebra $\cA$ by a finite group $G$ of $*$-automorphisms of $\cA$. It is shown that if a Positivstellensatz of type I resp. II is valid for the matrix algebra $M_n(\cA)$, then it is valid for the cross product algebra $\cA \times_\alpha G$ as well. Here the main technical tool is a conditional expectation from the matrix algebra $M_n(\cA)$ onto the algebra $\cA \times_\alpha G$. In Section \ref{sect_sos_com_matr} we give a new approach to Artin's theorem for matrix algebras over commutative $*$-algebras which is based on diagonalization of matrices by means of quasi-unitary matrices. In Section \ref{sect_sos_interval} we use conditional expectations to describe  matrix polynomials which are positive semidefinite on intervals $[a,b]$ and $[a,+\infty)$. In Section \ref{sect_sos_matr_field} we derive a Positivstellensatz of type I for the matrix algebra over formally real field equipped with the diagonal involution. Sections \ref{sect_sos_cyc} and \ref{sect_exam} deal with cyclic $*$-algebras. In Section \ref{sect_sos_cyc} we prove a Positivstellensatz of type IV for cyclic $*$-algebras. It states that positive elements belong to the quadratic module generated by (certain) finite products of commuting squares. In Section \ref{sect_exam} we elaborate a number of examples. The second example contains a positive element for which the Positivstellensatz of type II is not valid, but there is a Positivstellensatz with finite products of commuting squares in the numerator. Recently a question of Procesi and Schacher \cite{ps} was answered in the negative in \cite{klep}. Our third example provides another counterexample to this question. In the final Section \ref{openproblems} we list a number of open problems concerning non-commutative Positivstellens\"atze. 

\section{Preliminaries}\label{sect_prelim}

%The central notion is that of a $*$-algebra. 
%\medskip

Let $\bR$ be a formally real field and $\bC=\bR(\sqrt{-1})$ and let $\img=\sqrt{-1}$ be the imaginary unit. We define an involution on $\bC$ via $(l_1+\img l_2)^*=l_1-\img l_2,\ l_1,l_2\in\bR,$ and we denote by $\bR_+\subseteq\bR$ the set of positive elements of $\bR,$ or equivalently, of finite sums of squares in $\bR.$ For standard notions such as formally real fields, orderings, preorderings etc. we refer to the monographs \cite{pd}, \cite{mar} or \cite{schl}.
\begin{defn}\label{defn_star_alg}
An associative algebra $\cA$ over $\bR$ (resp. over $\bC$) is called a {\it $*$-algebra} if there is a map $x\mapsto x^*$ on $\cA$ called {\it involution} such that for $a,b\in\cA, \lambda\in\bR$ (resp. $\lambda\in\bC$):
\begin{enumerate}
 \item[\rm (i)] $(\lambda a+b)^*=\lambda a^*+b^*,$ (resp. $(\lambda a+b)^*=\lambda^*a^*+b^*$),
 \item[\rm (ii)] $(ab)^*=b^*a^*,$
 \item[\rm (iii)] $(a^*)^*=a.$
\end{enumerate}
\end{defn}

\noindent All $*$-algebras in this paper have an identity element denoted by $\Un{\cA}$ or simply by $\Un{}.$ 
%Simple computations show that $\Un{}^*=\Un{}.$ 

Let $\cV$ be a {\it unitary space}. That is, $\cV$ is a vector space over $\dC$ equipped with a scalar product $\langle\cdot,\cdot\rangle$ which is linear in the first variable and anti-linear in the second. Let $\cL(\cV)$ denote the space of all linear mappings of $\cV$ into itself.
\begin{defn}
An associative subalgebra $\cA\subseteq\cL(\cV)$ is called an $O^*$-algebra if for every $A\in\cA$ there exists  $B_A\in\cA$ such that $\langle Av,w\rangle=\langle v,B_Aw\rangle$ for all $v,w\in\cV.$ 
\end{defn}

Then the element $B_A$ is uniquely determined by $A$ and the $O^*$-algebra $\cA$ becomes a $*$-algebra with involution $A\to A^*:=B_A.$ Following \cite{s1}, $\cL^+(\cV)$ denotes the largest subalgebra of $\cL(\cV)$ which is an $O^*$-algebra.

\begin{defn} 
Let $\cA$ be a $*$-algebra over $\dC$ and $\cV$ be as above. A $*$-representation of $\cA$ on $\cV$ is a homomorphism $\pi:\cA\to\cL^+(\cV).$
\end{defn}

\noindent For a $*$-algebra $\cA$ let $\cA_h=\set{a\in\cA,\ a=a^*}$ denote the set of all self-adjoint elements of $\cA.$

\begin{defn}
Let $\cA$ be a $*$-algebra over $\bR$ or $\bC.$ A subset $\cC\subseteq\cA_h$ is a {\it quadratic} module if 
\begin{enumerate}
 \item[\rm (i)] $x+\lambda y\in\cC$ for  $x,y\in\cC,\ \lambda\in\bR_+,$
 \item[\rm (ii)] $z^*xz\in\c C$ for $x\in\cC,\ z\in\cA,$
 \item[\rm (iii)] $\Un{\cA}\in\cC.$
\end{enumerate}
\end{defn}

\noindent Let $\sum\cA^2$ be the smallest quadratic module in $\cA.$ It consists of elements of the form $x_1^*x_1^{}+\dots+ x_m^*x_m^{}$ which are called {\it sums of squares}.

\begin{defn}\label{defn_nc_pr}
A quadratic module $\cC\subseteq\cA_h$ is called a {\it non-commutative preordering} if
$c_1c_2\in\cC\ \mbox{for all}\ c_1,c_2\in\cC,\ c_1c_2=c_2c_1.$
\end{defn}

We denote by $\sum_{nc}\cA^2$  the smallest non-commutative preordering and call the elements of $\sum_{nc}\cA^2$  {\it non-commutative sums of squares}. If $\cA$ is a commutative $*$-algebra then $\sum\cA^2$ coincides with $\sum_{nc}\cA^2,$ but in general we have $\sum\cA^2\neq\sum_{nc}\cA^2.$

\begin{defn}
Let $a\in \cA_h.$ A subset $\cS_a\subseteq\cA_h$ is called a {\it denominator set} if
\begin{enumerate}
 \item[\rm (i)] $a\in\cS_a,$
 \item[\rm (ii)] if $b\in\cS_a$ and $x\in\cA$, then $x^*bx \in \cS_a,$
 \item[\rm (iii)] if $c\in\sum_{nc}\cA^2$ commutes with $b\in\cS_a$, then $cb \in \cS_a$.
\end{enumerate}
% $\cS_a$ is called a {\it weak denominator set} if in addition:
% \begin{enumerate}
%  \item[\rm (iv)] $a,b\in\cS_a$ implies $a+b\in\cS_a,$
% \end{enumerate}
\end{defn}

The preceding definitions are motivated by the following well-known  simple fact.
\begin{lemma}\label{lemma_prod_com}
Suppose $\cA$ be a (complex) $*$-algebra  of all {\it bounded} operators on a (complex) Hilbert space and $\cA_+$ be the set of positive operators in $\cA.$ Then we have $\sum_{nc}\cA^2\subseteq\cA_+$. If $a\in\cA_+$ then $\cS_a\subseteq\cA_+.$
\end{lemma}
\begin{proof}
It suffices to show that $ab\in\cA_+$ when $a,b\in\cA_+$ and $ab=ba$. By the functional calculus for bounded self-adjoint operators (see e.g. Theorem VII.1 in \cite{rs}), there exists a self-adjoint operator $b^{1/2}$ such that $(b^{1/2})^2=b$ and $ab^{1/2}=b^{1/2}a.$ Since $a\in\cA_+,$  we have
$$
\langle ab\varphi,\varphi\rangle=\langle b^{1/2}ab^{1/2}\varphi,\varphi\rangle=\langle ab^{1/2}\varphi,b^{1/2}\varphi\rangle\geq 0
$$
for $\varphi\in\cH$, that is $ab\in\cA_+.$
\end{proof}

%\begin{remk}
%The previous assertion does no hold in general if $a$ and $b$ are unbounded. 
%Namely, there exists a unitary space $\cV$ and $a,b\in L(\cV)$ such that for all %$\varphi\in\cV$ holds $\langle a\varphi,\varphi\rangle\geq 0,\langle %b\varphi,\varphi\rangle,\ ab\varphi=ba\varphi$ and $\langle %ab\varphi_0,\varphi_0\rangle<0$ for some $\varphi_0\in\cV.$
%\end{remk}

\begin{defn}\label{defn_cond_exp}
Let $\cA$ be a $*$-algebra over $\bR$ or $\bC.$  A linear map $p:\cA\to\cB$ is called a \textit{conditional expectation} of $\cA$ onto $\cB$ if
\begin{enumerate}
  \item[(i)] $p(a^*)=p(a)^*,\ p(b_1ab_2)=b_1p(a)b_2\ \mbox{for all}\ a\in\cA,\ b_1,b_2\in\cB,\ p(\mathbf{1}_\cA)=\mathbf{1}_\cB,$ and 
  \item[(ii)] $p(\sum\cA^2)\subseteq\sum\cA^2\cap\cB.$
\end{enumerate}

\noindent A linear map $p$ satisfying only condition (i) is called a $\cB$-\textit{bimodule projection} of $\cA$ onto $\cB.$ A conditional expectation $p$ will be called a \textit{strong conditional expectation} if
\begin{enumerate}
  \item[(ii')] $p(\sum\cA^2)\subseteq\sum\cB^2.$
\end{enumerate}
\end{defn}

\noindent Conditional expectations  for general  $*$-algebras have been introduced in \cite{ss}. In \cite{rief} and\cite{ss} they are used for the study of induced $*$-representations. In this paper  conditional expectations are important tools to prove Positivstellens\"atze for $*$-algebras.

%We conclude this section with definition of matrices over $*$-algebras and a simpe lemma, %which wil be used below. 
Let $\bL$ denote $\bR$ or $\bC$ and let $\cA$ be a $*$-algebra over $\bL$. We denote by $M_n(\bL)$ be $*$-algebra of matrices over $\bL$ considered with the standard linear base $\set{E_{ij}}_{i,j=\overline{1,n}}$ and with multiplication and involution defined by
\begin{gather}\label{eq_E_ij}
E_{ij}\cdot E_{kl}=\delta_{jk}\cdot E_{il},\ (\lambda E_{ij})^*:=\lambda^*E_{ji},\ \mbox{for}\ i,j,k,l\in\set{1,\dots,n},\ \lambda\in\bL.
\end{gather}

We define $M_n(\cA):=\cA\otimes_\bL M_n(\bL)$ with involution given by $(\sum_{ij}a_{ij}\otimes E_{ij})^*=\sum a_{ij}^*\otimes E_{ji}.$ We conclude this section with an elementary lemma.
\begin{lemma}\label{lemma_sos_rank1}
Let $\cA$ be a $*$-algebra. Each element of $\sum M_n(\cA)^2$ is a finite sum of "rank one" squares 
$$\sum_{ij}y_i^*y_j^{}\otimes E_{ij}=\left(\sum_{i}y_i\otimes E_{ki}\right)^*\left(\sum_{j}y_j\otimes E_{kj}\right),\ y_i\in\cA,\ i=1,\dots,n,$$
where $k\in\set{1,\dots,n}$ is fixed.
\end{lemma}
\begin{proof} The proof follows from the  computation
\begin{gather*}\left(\sum_{i,j=1}^n a_{ij}\otimes E_{ij}\right)^*\left(\sum_{i,j=1}^n a_{ij}\otimes E_{ij}\right)=\sum_{i=1}^n\left(\sum_{j=1}^n a_{ij}\otimes E_{ij}\right)^*\left(\sum_{j=1}^n a_{ij}\otimes E_{ij}\right)=\\ 
=\sum_{i=1}^n\left(\sum_{j=1}^n a_{ij}\otimes E_{kj}\right)^*\left(\sum_{j=1}^n a_{ij}\otimes E_{kj}\right).
\end{gather*}
\end{proof}

\section{Positivstellensatz of Type I for Path Algebras}\label{sect_sos_path} Many important $*$-algebras are quotients of path algebras of $*$-quivers, see e.g. \cite{cbh}. A Positivstellensatz for path algebras was recently proved by Popovych \cite{po}. The aim of this section to give an alternative proof of this result (Theorem \ref{pathpossatz}) by using a conditional expectation. 

Let $\Gamma=(\Gamma_0,\Gamma_1)$ be a $*$-quiver (an equivalent notion a $*$-double of a quiver was considered in \cite{cbh}, \cite{po} etc.). This means that $\Gamma$ is a directed multi-graph (i.e. multiple arrows between the vertices and knots are allowed) with a finite set of vertices $\Gamma_0=\set{e_1,\dots,e_n}$ and a finite set of arrows $\Gamma_1.$ For each $e_i,e_j\in\Gamma_0$ let $\Gamma(e_i,e_j)$ denote the set of arrows from $e_i$ to $e_j.$ For an arrow $b\in\Gamma(e_i,e_j)$ we denote by $o(b):=e_i$ and $t(b):=e_j$ the origin and the terminal vertex of $b$, respectively. 
For each arrow $b\in\Gamma(e_i,e_j)$ there exists a unique arrow $b^*\in\Gamma(e_j,e_i)$ and we have $(b^*)^*=b.$ For a knot $b\in\Gamma(e,e),\ e\in\Gamma_0$, we assume that $b^*\neq b.$ 

A path in $\Gamma$ is a finite sequence of arrows $b_1b_2\dots b_k$ such that $t(b_i)=o(b_{i+1}),\ i=1,\dots,k-1.$ We consider each vertex $e_i\in\Gamma_0$ as a path of the length zero. Let $\cB$ denote the union of the set of all paths in $\Gamma$ with $0.$ For two pathes $b_1b_2\dots b_m$ and $c_1c_2\dots c_l$ we define their product to be $b_1b_2\dots b_mc_1c_2\dots c_l$ if $t(b_m)=o(c_1)$ and $0$ otherwise. Then $\cB$ becomes a semigroup with respect to this multiplication. The semigroup algebra of this semigroup $\cB$ with involution determined by $b\mapsto b^*$ for $b\in\Gamma_0$ is called the {\it path algebra $\dC\Gamma$.}

As in \cite{po} we define an embedding $\epsilon:\dC\Gamma\to M_n(\cF),$ where $\cF$ is the free $*$-algebra with generator set $\Gamma_1$, by 
$$\epsilon(b):=b\otimes E_{ij}\ \mbox{for}\ b\in\Gamma(e_i,e_j)\ \mbox{and}\ \epsilon(e_i):=1\otimes E_{ii}.$$ 
%Using $\epsilon$ one proves the following Lemma (cf. %\cite{po}).
%\begin{lemma}
%Algebra $\dC\Gamma$ has a separating family of finite-dimensional $*$-representations.
%\end{lemma}
%\begin{proof}
%It is known that the free $*$-algebra $\cF$ has a separating family of %$*$-representations, see e.g. \cite{pu}. Thus, the same is true for $M_n(\cF)$ and %$\epsilon(\dC\Gamma).$
%\end{proof}
Our proof of Theorem \ref{pathpossatz} below uses the following slight generalization of Helton's theorem \cite{he}.
\begin{prop}\label{matrixfn}
Let $\cF_m=\dC\langle a_1^{},\dots,a_m^{},a_1^*,\dots,a_m^*\rangle$ be the free $*$-algebra with $m$ generators and let $X=X^*\in M_n(\cF_m).$ Then $\rho(X)\geq 0$ for every finite-dimensional $*$-representation $\rho$ of $M_n(\cF_m)$ if and only if $X\in\sum M_n(\cF_m)^2.$
\end{prop}
To prove this proposition we need the following technical result (see e.g. Lemma 2, \cite{s4}).
\begin{lemma}\label{closcone1}
Let $\cA$ be a unital $\ast$-algebra which has a faithful $\ast$-representation $\pi$ (that is, $\pi(a)=0$ implies that $a=0$) and is a union of a sequence of finite dimensional
subspaces $E_n$, $n \in \dN$. Assume that for each $n \in \dN$ there exists a number $k_n\in \dN$ such that the following is satisfied:
If $a \in \sum \cA^2$ is in $E_n$, then we can write $a$ as a finite sum $\sum_j ~a_j^\ast a_j$ such that all $a_j$ are in $E_{k_n}$.\\
Then the cone $\sum \cA^2$ is closed in $\cA$ with respect to the finest locally convex topology on $\cA$.
\end{lemma}
\noindent{\it Proof of Proposition \ref{matrixfn}.}
Since the if part is trivial, it suffices to prove the only if part. The main step of this proof is to show that the cone $\sum M_n(\cF_m)^2$ is closed in the finest locally convex topology on $M_n(\cF_m)$. For this we apply Lemma \ref{closcone1} to the $*$-algebra $\cA:=M_n(\cF_m)$.
The $*$-algebra $\cF_m$ has a faithful $*$-representation (see e.g. \cite{s4}), so has $M_n(\cF_m)$. Let $E_k$ be the vector space of generated by the elements $w\otimes E_{ij}$, $i,j=1,\dots,n,$ where $w$ runs over the monomials in $\cF_m$ of degree $\leq k.$ 
Let $A=\sum A_j^*A_j\in E_{2k}.$ Comparing the degrees of elements in the main diagonal of $A$, we conclude that each $A_j$ is in $E_k$. Hence the assumptions of Lemma \ref{closcone1} are fulfilled, so $\sum M_n(\cF_m)^2$ is closed.

Now we proceed almost verbatim as in the proofs of Propositions 4 and 5 in \cite{s4}. Assume to the contrary that $X\notin \sum M_n(\cF_m)^2$. By the separation theorem for convex sets there exists a linear functional $f$ on $M_n(\cF_m)$ such that $f$ is nonnegative on $\sum M_n(\cF_m)^2$ and $f(X)<0$. If $\pi_f$ denotes the representation of $M_n(\cF_M)$ obtained by the GNS-construction from $f$, there is a vector $\varphi $ of the representation space $V$ such that $f(A)=\langle \pi_f(A)\varphi,\varphi \rangle$ for all $A\in M_n(\cF_m)$. 
Let $P$ be the projection of $\cV$ onto
the finite-dimensional subspace $\pi(E_{2k})\varphi$. Since $a_1,\dots,a_m$ are generators of the free algebra $\cF_m$, there is a finite dimensional $\ast$-representation $\rho$ 
of $M_n(\cF_m)$ on $P \cV$ defined by $\rho(a_j)v=P\pi_f(a_j)v$, $v\in P\cV$, $j{=}1,\dots,m$. By construction we have $\pi_f(A)\varphi =\rho
(A) \varphi$ and hence $\langle\pi(B)\varphi,\varphi\rangle=\langle\rho(B)\varphi,\varphi\rangle$ for all $B\in E_{2k}.$ In particular, $f(x)= \langle\rho(X)\varphi,\varphi \rangle <0$ which contradicts the assumption. \hfill $\Box$
\medskip

Our next aim is to construct a strong conditional expectation from $M_n(\cF)$ onto $\dC\Gamma.$ For each $i,j=1,\dots,n$ we define a linear mapping $\gp_{ij}:\cF\to\cF$ as follows. For an element $b_1b_2\dots b_k\in\cF$, where $ b_1,b_2,\dots,b_k\in\Gamma_1$, put $\gp_{ij}(b_1b_2\dots b_k)=b_1b_2\dots b_k$ if $b_1b_2\dots b_k$ is a path in $\Gamma$ from $e_i$ to $e_j,$ that is, $$t(b_i)=o(b_{i+1}),\ i=1,\dots,k-1,\ o(b_1)=e_i,\ t(b_k)=e_j.$$ Otherwise we set $\gp_{ij}(b_1b_2\dots b_k)=0.$ Also we set $\gp_{ij}(\Un{}):=\Un{}.$ 
We will need the following auxiliary 
\begin{lemma}\label{lemma_gp_ij} For all $i,j=1,\dots,n $ and$\ x,y\in\cF$ we have:
\begin{enumerate}
  \item[(i)] $\gp_{ij}(x^*)=\gp_{ji}(x)^*,$
  \item[(ii)] $\gp_{ij}(x^*y)=\sum_{k=1}^n \gp_{ki}(x)^*\gp_{kj}(y).$
\end{enumerate}
\end{lemma}
\begin{proof}
Both equations follow directly from the definition of $\gp_{ij}.$
\end{proof}

For an arbitrary element $X=\sum_{ij}x_{ij}\otimes E_{ij},\ x_{ij}\in\cF$ of $M_n(\cF),$ we define $$\gP(X):=\sum_{ij}\gp_{ij}(x_{ij})\otimes E_{ij}.$$

\begin{prop}
The mapping $\gP$ is a strong conditional expectation from $M_n(\cF)$ onto $\dC\Gamma$.
\end{prop}
\begin{proof}
Some easy computations show that $\gP$ is a $\dC\Gamma$-bimodule projection. We prove the strong positivity property of $\gP.$ For let $Y\in \sum M_n(\cF)^2.$ By Lemma \ref{lemma_sos_rank1}, $Y$ is a finite sum of "rank one" squares, that is, we have $Y=\sum_{i,j}y_i^*y_j^{}\otimes E_{ij}$, where $y_i\in\cF.$ Using the definition of $\gP$ and Lemma \ref{lemma_gp_ij}, (ii) we compute
\begin{gather*}
\gP(\sum_{i,j=1}^n y_i^*y_j^{}\otimes E_{ij})=\sum_{i,j=1}^n \gp_{ij}(y_i^*y_j^{})\otimes E_{ij}=\sum_{i,j=1}^n\sum_{k=1}^n\gp_{ki}(y_i)^*\gp_{kj}(y_j^{})\otimes E_{ik}\cdot E_{kj}=\\
=\sum_{k=1}^n\sum_{i,j=1}^n\left(\gp_{ki}(y_i)\otimes E_{ki}\right)^*\left(\gp_{kj}(y_j)\otimes E_{kj}\right)=\\
=\sum_{k=1}^n\left(\sum_{i=1}^n\gp_{ki}(y_i)\otimes E_{ki}\right)^*\left(\sum_{j=1}^n\gp_{kj}(y_j)\otimes E_{kj}\right)\in\sum\dC\Gamma^2.
\end{gather*}
\end{proof}

Combining the preceding two propositions we obtain the following Positivstellensatz.
\begin{thm}\label{pathpossatz}
Let $X=X^*\in\dC\Gamma.$ Then $\pi(X)\geq 0$ for every finite-dimensional $*$-representation $\pi$ of the path algebra $\dC\Gamma,$ if and only if $X=\sum_{j=0}^k X_j^*X_j^{}$ for some elements $X_j\in\dC\Gamma.$
\end{thm}

\section{Artin's Theorem for Matrices over Noncomutative $*$-Algebras}\label{sect_nc_matr}

Throughout this section we suppose that $\cA$ is a unital $*$-algebra without zero divisors such that $\cA^\circ:=\cA\backslash\set{0}$ satisfies the left Ore condition (that is, given $a\in\cA,s\in\cA^\circ,$ there exist $b\in\cA,t\in\cA^\circ$ such that $ta=bs.$) 

We denote by $\cD_n(\cA)$ the diagonal matrices of $M_n(\cA),$ by $\cD_n(\cA)^\circ$ the diagonal matrices of $M_n(\cA)$ with non-zero entries on the diagonal, and by $\cL_n(\cA)$ the matrices $X=(x_{ij})\in M_n(\cA)$ such that $x_{ij}=0$ for $i\neq j$ and $x_{ii}\neq 0$ for all $i.$

\subsection{} In this first subsection we develop a general diagonalization procedure for hermitian matrices over $\cA.$ It might be of some interest in itself.

Let $a\in\cA,\gb=(b_1,\dots,b_n)\in M_{1,n}(\cA),\ \gb^*:=(b_1^*,\dots,b_n^*)^t\in M_{n,1}(\cA)$ and $C=(c_{ij})\in M_n(\cA)$ and consider a block matrix $A=A^*\in M_{n+1}(\cA)$ defined by 
\begin{gather}\label{eq_matr_abbC}
A=\left(
\begin{array}{ll}
a & \gb \\ 
\gb^* & C
\end{array}\right). 
\end{gather}

Assume that $a\neq 0.$ By the left Ore property of $\cA^\circ=\cA\backslash\set{0}$ each right fraction $b_i^*a^{-1}$ is a left fraction. All these left fractions can be brought to a common denominator. That is, there exist elements $s\in\cA^\circ$ and $f_1,\dots,f_n\in\cA$ such that $b_i^*a^{-1}=s^{-1}f_i$ or, equivalently $sb_i^*=f_ia$ for $i=1,\dots,n.$ Set $\gf=(f_1,\dots,f_n)^t\in M_{n,1}(\cA)$. Since $A=A^*$ and hence $a=a^*,$ we have
\begin{gather}\label{eq_ore}
s\gb^*=\gf a\ \mbox{and}\ a\gf^*=\gb s^*.
\end{gather}

Let $sC s^*:=(sc_{ij}s^*),\ \gf=(f_1,\dots,f_n)\in M_{n,1}(\cA)$ and $\gf a\gf^*:=(f_iaf_j)\in M_n(\cA)$ and put 
\begin{gather}\label{defD}
D=sC s^*-\gf a\gf^*.
\end{gather}

\begin{lemma}
Let $x\in\cA$ and $Y\in M_n(\cA).$ If (\ref{eq_ore}) holds, then 
\begin{gather}
\left(
\begin{array}{ll}
\ \ x & 0 \\ 
-Y\gf & Ys
\end{array}
\right)\left(
\begin{array}{ll}
a & \gb \\ 
\gb^* & C
\end{array}
\right)\left(
\begin{array}{ll}
\ \ x & 0 \\ 
-Y\gf & Ys
\end{array}
\right)^*=\left(
\begin{array}{ll}
xax^* & \ 0 \\ 
\ \ 0 & YDY^*
\end{array}
\right)
\end{gather}
\end{lemma}
\begin{proof}
We compute the matrix on the left-hand side by applying equation (\ref{eq_ore}) several times and obtain
\begin{gather*}
 \left(
\begin{array}{ll}
\ \ \ \ xa & \ \ \ \ x\gb \\ 
-Y\gf a+Ys\gb^* & -Y\gf\gb+YsC
\end{array}
\right)\left(
\begin{array}{ll}
\ \ x & 0 \\ 
-Y\gf & Ys
\end{array}
\right)^*=\left(
\begin{array}{ll}
xa & \ \ \ \ xb \\ 
0 & YsC-Y\gf\gb
\end{array}
\right)\left(
\begin{array}{ll}
x^* & -\gf^*Y^* \\ 
0 & \ \ s^*Y^*
\end{array}
\right)=\\
=\left(
\begin{array}{ll}
xax^* & -xa\gf^*Y^*+x\gb s^*Y^* \\ 
\ \ 0 & YsCs^*Y-Y\gf\gb s^*Y^*
\end{array}
\right)=\left(
\begin{array}{ll}
xax^* & \ \ \ \ \ \ \ \ \ \ 0 \\ 
\ \ 0 & YsCs^*Y^*-Y\gf a\gf^*Y^*
\end{array}
\right)=\left(
\begin{array}{ll}
xax^* & \ \ 0 \\ 
\ \ 0 & YDY^*
\end{array}
\right).
\end{gather*}

\end{proof}

Now we specialize the element $x\in\cA$ and the matrix $Y\in M_n(\cA).$ Suppose $x\neq 0.$ Applying once more the left Ore property we can write all right fractions $(Y\gf)_ix^{-1}\equiv\sum_iy_{ij}f_jx^{-1},\ i=1,\dots,n$, as left fractions with a common denominator, that is, there are elements $u\in\cA^\circ$ and $g_1,\dots,g_n\in\cA$ such that $(Y\gf)_ix^{-1}=u^{-1}g_i$ for $i=1,\dots,n.$ Setting $\gog=(g_1,\dots,g_n)^t,$ we have
\begin{gather}\label{eq_uYf=gx}
 uY\gf=\gog x.
\end{gather}

\begin{lemma}
If (\ref{eq_ore}) and (\ref{eq_uYf=gx}) are satisfied, then 
\begin{gather}
\left(
\begin{array}{ll}
x & \ 0 \\ 
0 & uYs
\end{array}
\right)
\left(
\begin{array}{ll}
a & \gb \\ 
\gb^* & C
\end{array}
\right)
\left(
\begin{array}{ll}
x & \ 0 \\ 
0 & uYs
\end{array}
\right)^*=\left(
\begin{array}{ll}
1 & 0 \\ 
\gog & uI_n
\end{array}
\right)
\left(
\begin{array}{ll}
xax^* & \ \ 0 \\ 
\ \ 0 & YDY^*
\end{array}
\right)
\left(
\begin{array}{ll}
1 & 0 \\ 
\gog & uI_n
\end{array}
\right)^*
\end{gather}
\end{lemma}
\begin{proof}
Let us denote by $L$ and $R$ the matrices on the left and right hand-sides, respectively. We compute the right-hande side and obtain
\begin{gather*}
R=\left(
\begin{array}{ll}
xax^* & \ \ \ 0 \\ 
\gog xax^* & uYDY^*
\end{array}
\right)
\left(
\begin{array}{ll}
1 & \ \gog^* \\ 
0 & u^*I_n
\end{array}
\right)
=\left(
\begin{array}{ll}
xax^* & \ \ \ \ \ xax^*\gog^* \\ 
\gog xax^* & \gog xax^*\gog^*+uYDY^*u^*
\end{array}
\right)
\end{gather*}
By (\ref{defD}) and (\ref{eq_uYf=gx}), we compute
\begin{gather*}
\gog xax^*\gog^*+uYDY^*u^*=uY\gf a(uY\gf)^*+uYDY^*u^*=uY(\gf a\gf^*+D)Y^*u^*=uYsCs^*Y^*u.
\end{gather*}
Therefore, by (\ref{eq_uYf=gx}) and (\ref{eq_ore}), we continue and derive 
\begin{gather*} 
R=\left(
\begin{array}{ll}
xax^* & xa\gf^*Y^*u^* \\ 
uY\gf ax & uYsCs^*Y^*u^*
\end{array}
\right)=\left(
\begin{array}{ll}
xax^* & x\gb(uYs)^* \\ 
uYs\gb^*x^* & uYsC(uYs)^*
\end{array}
\right)
=L
\end{gather*}
which proves the assertion of the Lemma.
\end{proof}

\noindent\textbf{Remark.} Retaining the preceding notations, we have 
\begin{gather}
 \left(
\begin{array}{ll}
1 & 0 \\ 
\gog & uI_n
\end{array}
\right)
\left(
\begin{array}{ll}
x & 0 \\ 
-Y\gf & Ys
\end{array}
\right)
=\left(
\begin{array}{ll}
x & \ 0 \\ 
0 & uYs
\end{array}
\right).
\end{gather}

\begin{lemma}
 Let $Y\in \cL_n(\cA).$ Then there is another matrix $T\in\cL_n(\cA)$ such that $TY\in\cD_n(\cA)^\circ.$
\end{lemma}
\begin{proof}
 We proceed by induction on $n.$ Suppose that the assertion is proved for $n$ and let $Y\in\cL_{n+1}(\cA).$ We write 
$$
\left(
\begin{array}{ll}
y_0 & 0 \\ 
\gz & Y_n 
\end{array}
\right)
$$
with $y_0\in\cA_0^\circ, \gz\in M_{n,1}(\cA)$ and $Y_0\in\cL_n(\cA).$ By induction hypothesis there is a matrix $T_n\in\cL_n(\cA)$ such that $T_nY_n\in\cL_n(\cA).$ Note that $y_0\neq 0$ by the definition of $\cL_k(\cA).$ By the Ore property there exist elements $t_i,s_i\in\cA^\circ$ such that $(T_n\gz )_iy_0^{-1}=s_i^{-1}t_i,\ i=1,\dots,n.$ Let $S$ be the diagonal matrix with entries $s_i$ and $t$ the row with entries $t_i.$ Since then $ty_0=ST_n \gz,$
$$
T=\left(
\begin{array}{ll}
\ 1 & \ 0 \\ 
-\gz & ST_n
\end{array}
\right)
$$
has the desired property.
\end{proof}

\subsection{}
From now on we suppose that $\cA$ is an $O^*$-algebra on a unitary space $\cV.$ Then $M_n(\cA)$ is an $O^*$-algebra acting on $\cV_n=\cV\oplus\dots\oplus\cV$ ($n$ times).% by $(Av)_i\sum_ja_{ij}v_j$
If $A$ is an element of $\cA$ resp. $M_n(\cA)$ and $\cE$ is a linear subspace of $\cV$ resp. $\cV_n$, we shall write $A\geq 0$ on $\cE$ when $\langle A\varphi,\varphi \rangle \geq 0$ for all $\varphi \in \cE$. Define
\begin{align}\label{defaplusmnaplus}
\cA_+=\{a \in \cA: a\geq 0 ~~{\rm on}~~\cV\},~~~M_n(\cA)_+=\{A\in M_n(\cA): A\geq 0 ~~{\rm on}~~\cV_n\}.
\end{align}

Now let $A$ be a matrix given by (\ref{eq_matr_abbC}) and retain the above notation.
Let $\cE$ and $\cF_n$ be linear subspaces of $\cV$ and $\cV_n,$ respectively.

\begin{lemma}\label{lemma_A>0_on_E_Fn}
 $A\geq 0$ on $(\cE,\cF_n)$ if and only if $a\geq 0$ on $\cE,\ C\geq 0$ on $\cF_n$ and 
\begin{gather}\label{eq_cauchy_ineq}
 |\langle \gb\varphi,\varphi_1\rangle|^2\leq\langle a\varphi_,\varphi_1\rangle\langle C\varphi,\varphi\rangle\ \mbox{for}\ \varphi_1\in\cE,\ \varphi\in\cF_n.
\end{gather}
\end{lemma}

\begin{proof}
 Let $\alpha$ and $\beta$ be complex numbers and put $\psi_{\alpha,\beta}:=(\alpha\varphi_1,\beta\varphi).$ Then we compute 
\begin{gather}\label{eq_A_psi_psi}
\langle A\psi_{\alpha,\beta},\psi_{\alpha,\beta}\rangle=
\alpha\overline{\alpha}\langle a\varphi_1,\varphi_1\rangle+\alpha\overline{\beta}\langle\gb^*\varphi,\varphi_1\rangle+\overline{\alpha}\beta\langle \gb\varphi,\varphi_1\rangle+\beta\overline{\beta}\langle C\varphi,\varphi\rangle.
\end{gather}
Clearly, $A\geq 0$ on $(\cE,\cF_n)$ if and only if $a\geq 0$ on $\cE,\ C\geq 0$ on $\cF_n$ and $\langle A\psi_{\alpha,\beta},\psi_{\alpha,\beta}\rangle\geq 0$ for all $\varphi_1\in\cE,\varphi\in\cF_n$ and $\alpha,\beta\in\dC.$ Since the numbers $\alpha,\beta\in\dC$ in equation (\ref{eq_A_psi_psi}) are arbitrary, 
%$\alpha,\beta\in\dC,$ 
it follows that the latter is equivalent to the inequality (\ref{eq_cauchy_ineq}) as stated in the Lemma.
\end{proof}

\begin{cor}\label{cor_1}
 If $A\geq 0$ on $\cV_{n+1}$ and $a=0,$ then $\gb=0.$
\end{cor}

\begin{lemma}\label{lemma_abbC>0}
\begin{enumerate}
 \item[(i)] If $A\geq 0$ on $\cE_{n+1},$ then $a\geq 0$ on $\cE,\ C\geq 0$ on $\cE_n$ and $D\geq 0$ on $\cE_n.$
 \item[(ii)] If $a\geq 0$ on $\cE$ and $D\geq 0$ on $\cE,$ then $A\geq 0$ on $(\gf^*\cE,s^*\cE).$
\end{enumerate}
\end{lemma}
\begin{proof}(i): Suppose that $A\geq 0.$ Then $a\geq 0$ on $\cE$ and $C\geq 0$ on $\cE_n$ by Lemma \ref{lemma_A>0_on_E_Fn}. Let $\varphi\in\cE_n.$ Using the identity $a\gf^*=bs^*$ by (\ref{eq_ore}) and inequality (\ref{eq_cauchy_ineq}) we conclude that 
\begin{gather*}
 |\langle \gf a\gf^*\varphi,\varphi\rangle|^2=|\langle \gf \gb s^*\varphi,\varphi\rangle|^2=|\langle \gb s^*\varphi,\gf^*\varphi\rangle|^2\leq\langle a\gf^*\varphi,\gf^*\varphi\rangle
\langle Cs^*\varphi,s^*\varphi\rangle=\langle \gf a\gf^*\varphi,\varphi\rangle\langle sCs^*\varphi,\varphi\rangle.
\end{gather*}

If $\langle \gf a\gf^*\varphi,\varphi\rangle\neq 0,$ then we have $\langle \gf a\gf^*\varphi,\varphi\rangle\leq \langle sCs^*\varphi,\varphi\rangle$ and so $\langle D\varphi,\varphi\rangle\geq 0.$ If $\langle \gf a\gf^*\varphi,\varphi\rangle=0,$ then $\langle D\varphi,\varphi\rangle=\langle Cs^*\varphi,s^*\varphi\rangle\geq 0,$ because $C\geq 0.$ Thus, $D\geq 0$ on $\cE_n.$

(ii): Let $\varphi,\psi\in\cE_n.$ Since $a\geq 0$ on $\cE,$ we have $\langle D\psi,\psi\rangle=\langle Cs^*\psi,s^*\psi\rangle-\langle \gf a\gf^* \psi,\psi\rangle\geq 0$ and hence $\langle Cs^*\psi,s^*\psi\rangle\geq\langle a\gf^*\psi,\gf^*\psi\rangle=\langle \gb s^*\psi,\gf^*\psi\rangle\geq 0,$ so we obtain $C\geq 0$ on $s^*\cE_n$ and 
\begin{gather*}
|\langle \gb s^*\varphi,\gf^*\psi\rangle|^2=|\langle a\gf^*\varphi,\gf^*\psi\rangle|^2\leq\langle a\gf^*\varphi,\gf^*\varphi\rangle\langle a\gf^*\psi,\gf^*\psi\rangle\leq\langle a\gf^*\varphi,\gf^*\varphi\rangle\langle Cs^*\psi,s^*\psi\rangle.
\end{gather*}

Therefore, $A\geq 0$ on $(\gf^*\cE,s^*\cE)$ by Lemma \ref{lemma_A>0_on_E_Fn}.
\end{proof}

\begin{prop}\label{prop_diagonalisation}
For each matrix $A\in M_n(\cA)_+$ there exist matrices $X_+,X_-\in\cL_n(\cA)$ such that $X_+^{}AX_+^*\in\cD_n(\cA)_+$ and $X_-X_+\in\cD_n(\cA)^\circ.$
\end{prop}
\begin{proof}
Let $z_1,\dots,z_n$ be given elements of $\cA^\circ.$ In view of the subsequent application given below we prove the stronger assertion that the diagonal matrix $\cD=X_+AX_+^*$ can be chosen the form $d_i=z_ia_iz_i^*$ for some $a_i\in\cA_+.$

The proof is given by induction on $n.$ Obviously, the assertion is true for $n=1.$ Let $A\in M_{n+1}(\cA)_+.$ We write $A$ in the form (\ref{eq_matr_abbC}). Since $A\geq 0,$ we have $C\geq 0$ by Lemma \ref{lemma_abbC>0}, so the induction hypothesis applies to the matrix $C.$

If $a=0,$ then $\gb=0$ by Corollary \ref{cor_1} and it suffices to enlarge the corresponding matrices for $C$ by putting $1$ in the left upper corner and $0$ elsewhere.

From now on suppose that $a\neq 0.$ By the induction hypothesis, there are matrices $Y_+,Y_-\in\cL_n(\cA)$ for which $Y_+Y_-\in\cD_n(\cA)^\circ$ and $D_n:=Y_+CY_+^*\in\cD_n(\cA)_+$ has diagonal entries $d_i=z_ia_iz_i^*$ with $a_i\in\cA_+,\ i=2,\dots,n+1.$ We apply Lemma \ref{lemma_abbC>0} with $x=z_1$ and $Y=Y_+.$ Putting 
$$
X_+=
\left(
\begin{array}{ll}
\ \ z_1 & 0 \\ 
-Y_+\gf & Y_+s
\end{array}
\right)
$$
we therefore have $X_+\in\cL_{n+1}(\cA)$ and $X_+AX_+\in\cD_{n+1}(\cA)_+$ has the diagonal $z_1az_1^*,z_ia_iz_i^*$ for $i=2,\dots,n+1.$ Note that $a\in\cA_+.$

From Lemma \ref{lemma_A>0_on_E_Fn} there is a matrix $T\in\cL_n(\cA)$ such that $T\cdot(uY_-s)\in\cD_n(\cA)^\circ.$ Set 
$$
X_-=\left(
\begin{array}{ll}
1 & 0 \\ 
0 & T
\end{array}
\right)\left(
\begin{array}{ll}
1 & 0 \\ 
\gog & uI_n
\end{array}
\right).
$$
From equation (\ref{eq_uYf=gx}) it follows that 
\begin{gather*}
 X_-X_+=\left(
\begin{array}{ll}
1 & 0 \\ 
0 & T
\end{array}
\right)
\left(
\begin{array}{ll}
1 & 0 \\ 
\gog & uI_n
\end{array}
\right)
\left(
\begin{array}{ll}
\ \ z_1 & 0 \\ 
-Y_+\gf & Y_+s
\end{array}
\right)
=\left(
\begin{array}{ll}
z_1 & 0 \\ 
0 & TuYs
\end{array}
\right)
\in\cD_n(\cA)^\circ.
\end{gather*}

\end{proof}

\subsection{ } Let $A\in M_n(\cA)_+.$ Then, by Proposition \ref{prop_diagonalisation} there are diagonal matrices $D_0\in\cD_n(\cA)_+,\ D\in\cD_n(\cA)^\circ$ and matrices $X_+,X_-\in\cL_n(\cA)$ such that 
\begin{gather}\label{eq_D_+=X_+AX_+*}
 D_0=X_+AX_+^*\ \mbox{and}\ DAD^*=X_-D_0X_-^*.
\end{gather}

We shall use this result to show that a Positivstellensatz of type $II$ holds for the matrices over $\cA$ provided that it holds for $\cA$ itself. More precicely, we have the following Positivstellensatz. Recall that $\cA_+$ and $M_n(\cA)_+$ have been defined by (\ref{defaplusmnaplus}).
\begin{thm}\label{posncmatrices}
Let $\cA$ be an $O^*$-algebra. Suppose that $\cA$ has no zero divisors and $\cA$ satisfies the left Ore condition. Assume that for each element $a\in\cA_+$ there exists $z\in\cA^\circ$ such that $zaz^*\in\cA^2.$

Then for each matrix $A\in M_n(\cA)_+,\ n\in\dN,$ there are matrices $D_0\in\sum\cD_n(\cA)^2,D\in\cD_n(\cA)^\circ$ and $X_+,X_-\in\cL_n(\cA)$ such that $X_+AX_+^*=D_0$ and $DAD^*=X_-D_0X_-^*\in \sum M_n(\cA)^2.$
\end{thm}
\begin{proof}
Let $z_1,\dots,z_n$ be fixed elements of $\cA_0.$ In the above proof of the Proposition \ref{prop_diagonalisation} it was shown that there exist elements $a_1,\dots,a_n$ such that $X_+AX_+^*\in\cD_n(\cA)_+$ has the diagonal entries $z_iaz_i^*.$ By the assumption we can choose $z_i\in\cA$ such that $z_ia_iz_i^*\in\sum\cA^2.$ Then the assertion follows from the proof of Proposition \ref{prop_diagonalisation}, see also (\ref{eq_D_+=X_+AX_+*}).
\end{proof}
\textbf{Remark.} In Subsection 4.3 of \cite{s2} a related result was obtained for matrices over the {\it commutative} polynomial algebra $\dR[x_1,\dots,x_d]$. In this result we had matrices $X_+,X_-\in\cL_n(\cA)$ for which both products $X_+X_-$ and $X_-X_+$ are central. In the above theorem for the {\it noncommutative} $O^*$-algebras $\cA$ we have only the weaker assertion stating that $X_-X_+\in\cD_n(\cA)^\circ.$
\section{Positivstellens\"atze for Crossed Product Algebras}\label{sect_sos_cross}
Let $\cA$ be a unital $*$-algebra and let $G$ be a finite group of $*$-automorphisms of $\cA.$ Let $\alpha_g\in\Aut G$ denote the $*$-automorphism corresponding to $g\in G.$ In this section we will show how Positivstellens\"atze for the matrix algebra $M_n(\cA),\ n=|G|,$ can be used to derive Positivstellens\"atze for the crossed product algebra $\cA\times_\alpha G$ of $\cA$ with $G.$

First let us recall the defintion of the crossed product $*$-algebra $\cA\times_{\alpha}G$. As a linear space it is the tensor product $\cA\otimes\dC[G]$ or equivalently the vector space of $\cA$-valued functions on $G$ with finite support. Product and involution on $\cA$ are determined by 
\begin{align*}(a\otimes g)(b\otimes h)=a\alpha_g(b)\otimes gh~~{\rm and}~~(a\otimes g)^*=\alpha_{g^{-1}}
(a^*)\otimes g^{-1},
\end{align*}
 respectively. If we identify $b$ with $b\otimes e$ and $g$ with $1 \otimes g$, then the $*$-algebra $\cA\times_{\alpha}G$ can be considered as the universal $*$-algebra
generated by the two $*$-subalgebras $\cA$ and $\dC[G]$ with cross commutation relations $gb=\alpha_g(b)g$ for $b\in \cA$ and $g\in G$.

\subsection{} Our first aim is to construct an embedding $\cA\times_\alpha G\hookrightarrow \cA\otimes M_n(\dC).$ Define the linear mapping $\epsilon$ from $\cA\times_\alpha G$ to $M_n(\cA)$ as follows:
\begin{gather}\label{eq_emb_AG_MnA}
\epsilon:a\otimes g \mapsto \sum_{h\in G}\alpha_h(a)\otimes E_{h,hg}.
\end{gather}

\begin{lemma}\label{lemma_emb_AG_MnA}
The mapping $\epsilon$ is an injective $*$-homomorphism of $\cA\times_\alpha G$ to $M_n(\cA)$. 
\end{lemma}
\begin{proof}
Take $a\otimes g,\ b\otimes k\in\cA\times_\alpha G.$ Then we have
\begin{gather*}
\epsilon(a\otimes g)\epsilon(b\otimes k)=\left(\sum_{h\in G}\alpha_h(a)\otimes E_{h,hg}\right)\left(\sum_{l\in G}\alpha_l(b)\otimes E_{l,lk}\right)=\\
=\sum_{h,l\in G}\delta_{hg,l}\cdot\alpha_h(a)\alpha_l(b)\otimes E_{h,lk}=\sum_{h\in G}\alpha_h(a)\alpha_{hg}(b)\otimes E_{h,hgk}=\sum_{h\in G}\alpha_h(a\alpha_g(b))\otimes E_{h,hgk}=\\
=\epsilon(a\alpha_g(b)\otimes gk)=\epsilon((a\otimes g)(b\otimes k)).
\end{gather*}
Analogously one checks that $(\epsilon(a\otimes g))^*=\epsilon((a\otimes g)^*).$ Thus $\epsilon$ is a $*$-homomorphism. It is easily seen that $\epsilon$ is injective.
\end{proof}

From now on we consider $\cA\times_\alpha G$ as a $*$-subalgebra of $M_n(\cA)$ via the embedding $\epsilon.$

\subsection{} Next we define a projection $\gP$ from $M_n(\cA)$ onto $\cA\times_\alpha G.$ For every $g\in G$ let $\beta_g$ denote the linear mapping of $M_n(\cA)$ onto itself defined by 
\begin{gather}\label{eq_act_G_MnA}
\beta_g: a\otimes E_{m,k}\mapsto\alpha_g(a)\otimes E_{gm,gk},\ a\in\cA,\ m,k,g\in G.
\end{gather}

\begin{lemma}
The map $g\mapsto\beta_g$ 
%defined by (\ref{eq_act_G_MnA}) 
is a well-defined action of $G$ on $M_n(\cA)$ by $*$-automorphisms.
\end{lemma}
\begin{proof}
The proof is given by straightforward computations. We omit the detail.% show that $\beta_g((a\otimes E_{hl})(c\otimes E_{mk}))=\beta_g(a\otimes E_{hl})\beta_g(c\otimes E_{mk})$.
\end{proof}

Define $\gP$ as the average over the action $g\mapsto\beta_g,$ that is,
\begin{gather}
\gP(a\otimes E_{m,k}):=\frac{1}{n}\sum_{g\in G}\beta_g(a\otimes E_{m,k})=\frac{1}{n}\sum_{g\in G}\alpha_g(a)\otimes E_{gm,gk}.
\end{gather}

\begin{prop}\label{thm_str_c_exp}
The mapping $\gP$ is a faithful strong conditional expectation from $M_n(\cA)$ onto $\cA\times_\alpha G.$
\end{prop}
\begin{proof}
It follows from the formulas (\ref{eq_emb_AG_MnA}) and (\ref{eq_act_G_MnA}) that $\cA\times_\alpha G$ is the stable $*$-subalgebra under the action of $G$ on $M_n(\cA).$ Being an average over an action of a finite group $\gP$ is a faithful conditional expectation by \cite{ss}, Proposition 5.

We prove the strong positivity property of $\gP.$ Take an element $\sum_{m,k\in G}c_{m,k}\otimes E_{m,k}\in M_n(\cA).$ Then we have
\begin{gather}
\nonumber \gP\left(\left(\sum_{m,k\in G}c_{m,k}\otimes E_{m,k}\right)^*\left(\sum_{m,k\in G}c_{m,k}\otimes E_{m,k}\right)\right)=\gP\left(\sum_{m,k,l\in G}c_{m,k}^*c_{m,l}^{}\otimes E_{k,l}\right)=\\
\label{eq_long}=\frac{1}n\sum_{g\in G}\beta_g\left(\sum_{m,k,l\in G}c_{m,k}^*c_{m,l}^{}\otimes E_{k,l}\right)=\frac{1}n\sum_{g\in G}\sum_{m,k,l\in G}\alpha_g(c_{m,k}^*c_{m,l}^{})\otimes E_{gk,gl}=\\
\nonumber=\frac{1}n\sum_{m\in G}\left(\sum_{g,k\in G}\alpha_g(c_{m,k})\otimes E_{gm,gk}\right)^*\left(\sum_{g,l\in G}\alpha_g(c_{m,l})\otimes E_{gm,gl}\right)
\end{gather}
The embedding formula (\ref{eq_emb_AG_MnA}) implies that
\begin{gather*}
\sum_{g,k\in G}\alpha_g(c_{m,k})\otimes E_{gm,gk}=\sum_{k\in G}\sum_{g\in G}\alpha_{gm}(\alpha_{m^{-1}}(c_{m,k}))\otimes E_{gm,gm(m^{-1}k)}=\\
=\sum_{k\in G}\sum_{g\in G}\alpha_g(\alpha_{m^{-1}}(c_{m,k}))\otimes E_{g,gm^{-1}k}=\sum_{k\in G}\alpha_{m^{-1}}(c_{m,k})\otimes m^{-1}k.
\end{gather*}
Using the latter equation we proceed in (\ref{eq_long}) and derive
\begin{gather*}
=\frac{1}n\sum_{m\in G}\left(\sum_{k\in G}\alpha_{m^{-1}}(c_{m,k})\otimes m^{-1}k\right)^*\left(\sum_{k\in G}\alpha_{m^{-1}}(c_{m,k})\otimes m^{-1}k\right)\in\sum(\cA\times_\alpha G)^2.
\end{gather*}
\end{proof}

\subsection{} As in Section \ref{sect_nc_matr}, we suppose that $\cA$ is an $O^*$-algebra acting on a unitary space $\cV.$ 

Let us recall the definition of the regular covariant representation of $\cA\times_\alpha G$. The representation space $\cV^{|G|}$ is a direct sum $\oplus_{g\in G}\cV$ of $|G|$ copies of $\cV.$ For $e\in\cV$ and $\ k\in G$ we denote by $e_k$ the element of $\oplus_{g\in G}\cV$ which has $e$ at the place $k$ and is $0$ otherwise. Let $g \in G$ and $a\in \cA$. We define linear mappings $\rho(g)$, $\pi(a)$ and $\pi_{creg}(a\otimes g)$ on the vector space $\cV^{|G|}$ by 
$$
\rho(g)e_k:=e_{kg^{-1}},\ \ \pi(a)e_k:=(\alpha_k(a)e)_k,\ \ \pi_{creg}(a\otimes g)e_k=\pi(a)\rho(g)e_k,\ \ e\in\cV,\ k\in G.
$$
Some simple computations show that $\pi_{creg}$ is a well-defined $*$-representation of the $*$-algebra $\cA\times_\alpha G$ on $\cV^{|G|}$ and that, 
$$\pi(\alpha_g(a))=\rho(g)\pi(a)\rho(g)^*,\ a\in\cA,\ g\in G.$$
The $*$-representation $\pi_{creg}$ (or likewise the triple $(\pi,\rho,\cV^{|G|})$) are called the \textit{regular covariant $*$-representation} of the crossed product algebra $\cA\times_\alpha G$. Define
$$(\cA\times_\alpha G)_+=\{x=x^* \in \cA\times_\alpha G: \pi_{creg}(x)\geq 0\}.
$$
Using the embedding formula (\ref{eq_emb_AG_MnA}) one easily verifies that the action of $\cA\times_\alpha G\subseteq M_n(\cA)$ coincides with $\pi_{creg},$ thus proving the
\begin{lemma}\label{prop_pos_cr_pr}
 $(\cA\times_\alpha G)_+\subseteq M_n(\cA)_+,$ where $M_n(\cA)_+$ is defined by (\ref{defaplusmnaplus}). 
\end{lemma}

% \begin{remk}\label{remk_1}
% Let $\cA$ be an $O^*$-algebra on a domain $\cV$ and assume that a crossed product $\cA\times_\alpha G$ with a finite group $G$ is defined. Then natural projection $p:\cA\times_\alpha G\to\cA\otimes e$ is a strong conditional expectation, see \cite{ss}. If we induce the identity representation of the subalgebra $\cA\otimes e\simeq\cA,$ its matrix representation corresponds to the embedding defined by (\ref{eq_emb_AG_MnA}). This representation of $\cA$ is sometimes called Heisenberg representation, see \cite{s3}. Defining positivity via this representation we get $(\cA\times_\alpha G)_+\subseteq M_n(\cA)_+.$
% \end{remk}

\begin{thm}
If the Positivstellensatz of type I is valid for $M_n(\cA),$ then it holds also for $\cA\times_\alpha G.$
\end{thm}
\begin{proof} Let $x=x^*\in(\cA\times_\alpha G)_+$. Then we have $x\in M_n(\cA)_+$ by the preceding lemma and hence $x\in\sum M_n(\cA)^2$ by the Positivstellensatz for $M_n(\cA)_+$. Since $\gP$ is a strong conditional expectation by Proposition \ref{thm_str_c_exp}, we obtain $\gP(x)=x\in\sum(\cA\times_\alpha G)^2.$
\end{proof}

\begin{thm}
If the Positivstellensatz of type III holds for $M_n(\cA),$ it holds for the cross product algebra $\cA\times_\alpha G$ as well.
\end{thm}
\begin{proof} Suppose that $x=x^*\in(\cA\times_\alpha G)_+.$ Then $x\in M_n(\cA)_+.$ Since the Positivstellensatz of type III holds for $M_n(\cA)$, there exists an element $y\in M_n(\cA)^\circ\cap\sum M_n(\cA)^2$ such that $xy=yx$ and $xy\in\sum M_n(\cA)^2$. Because $\gP$ is a strong conditional expectation by Proposition \ref{thm_str_c_exp}, we have $x\gP(y)=\gP(xy)=\gP(yx)=\gP(y)x\in\sum (\cA\times_\alpha G)^2$ and $\gP(y)\in\sum(\cA\times_\alpha G)^2$.

It is left to show that $\gP(y)\in(\cA\times_\alpha G)^\circ.$ Moreover, we show that $\gP(y)\in M_n(\cA)^\circ.$ Indeed, by assumption $y=\sum_i y_i^*y_i^{}.$ Let $z\in M_n(\cA)$ be such that $\gP(y)z=0.$ Then using defintion of $\gP$ we have also 
$$
0=z^*\gP(y)z=\frac{1}n\sum_i\sum_{g\in G}z^*\beta_g(y_i)^*\beta_g(y_i)z=\frac{1}n\sum_i\sum_{g\in G}(\beta_g(y_i)z)^*(\beta_g(y_i)z).
$$
Since $M_n(\cA)$ as well as $\cA$ is an $O^*$-algebra, all $\beta_g(y_i)z$ are $0,$ hence $yz=\sum_i y_i^*(y_i^{}z)=0.$ Since $y$ is not a zero divizor we get $z=0.$
\end{proof}

\noindent Now we turn to the Positivstellensatz of type II for $\cA\times_\alpha G.$ 
\begin{thm}
Let $\cA$ be an $O^*$-algebra without zero divisors such that $\cA^\circ$ is a left Ore set. Assume that the Positivstellensatz of type II holds for $\cA.$ Then the Positivstellensatz of type II holds also for the crossed-product algebra $\cA\times_\alpha G$.
\end{thm}
\begin{proof}
%Let $|G|=n,\ G=\set{g_1,\dots,g_n}$ and, as before, consider $\cA\times_\alpha G$ as a %$*$-subalgebra of $M_n(\cA)$ via the embedding (\ref{eq_emb_AG_MnA}). 
Suppose that $X\in(\cA\times_\alpha G)_+$. Then $X\in M_n(\cA)_+$ by Lemma \ref{prop_pos_cr_pr}. Therefore, by Theorem \ref{posncmatrices} there exists a diagonal matrix $Y\in \cD_n(\cA)^\circ$ such that $Y^*XY\in\sum M_n(\cA)^2$. Let $Y=\sum_{i=1}^n y_i^{}\otimes E_{g_i,g_i},$ where $y_1,y_2,\dots,y_n\in\cA^\circ.$ Because $\alpha_{g_i^{-1}}(y_i)\in\cA^\circ,$ it follows from a repeated application of the Ore condition that there exist $z_1,z_2,\dots,z_n\in\cA^\circ$ such that 
\begin{gather}\label{eq_YZ_in_crpr}
\alpha_{g_1^{-1}}(y_1)z_1=\dots=\alpha_{g_n^{-1}}(y_n)z_n\in\cA^\circ.
\end{gather}
Put $Z=\sum_{i=1}^n\alpha_{g_i}(z_i)\otimes E_{g_i,g_i}.$ We claim that $YZ\in(\cA\times_\alpha G)^\circ.$ It is easily seen that $YZ\in M_n(\cA)^\circ.$ We check that $YZ$ is invariant under $\beta_g,\ g\in G,$ which implies $YZ\in\cA\times_\alpha G.$ Indeed, for a fixed $g\in G$ we compute using (\ref{eq_YZ_in_crpr})% and definition of $\beta_g$
\begin{gather*}
\beta_g(YZ)=\beta_g\left(\sum_{i=1}^ny_i\alpha_{g_i}(z_i)\otimes E_{g_i,g_i}\right)=\sum_{i=1}^n\alpha_g(y_i\alpha_{g_i}(z_i))\otimes E_{gg_i,gg_i}=\\
=\sum_{g_i\in G}\alpha_{gg_i}(\alpha_{g_i^{-1}}(y_i)z_i)\otimes E_{gg_i,gg_i}=\sum_{g_j\in G}\alpha_{g_j}(\alpha_{g_j^{-1}}(y_j)z_j)\otimes E_{g_j,g_j}=\sum_{j=1}^ny_j\alpha_{g_j}(z_j)\otimes E_{g_j,g_j}=YZ.
\end{gather*}
Since $(YZ)^*XYZ =Z^*Y^*XYZ\in\sum M_n(\cA)^2$ and $X, YZ\in(\cA\times_\alpha G)^\circ,$ we obtain $\gP((YZ)^*XYZ)=(YZ)^*XYZ\in\sum(\cA\times_\alpha G)^2.$
\end{proof}

\begin{exam}
Let us describe $\epsilon$ more explicitly in the case when $G=\set{0,1,\dots,n-1}$ is a cyclic group of order $n$ generated by a $*$-automorphism $\sigma,\ \sigma^n=\Id.$ Then $\epsilon$ maps the element $\sum_{k=0}^{n-1}a_k\otimes k$ onto the matrix
\begin{gather*}
\left(\begin{array}{llll}
a_0 & a_1 & \dots & a_{n-1} \\ 
\sigma(a_{n-1}) & \sigma(a_0) & \dots & \sigma(a_{n-2})\\ 
\vdots & \vdots & \ddots & \vdots \\ 
\sigma^{(n-1)}(a_1) & \sigma^{(n-1)}(a_2) & \dots & \sigma^{(n-1)}(a_0)
\end{array}\right)
\end{gather*}
The algebra of matrices of this form might be of interest in itself even if $\sigma$ is the identity automorphism.

\end{exam}

\section{Artins Theorem for Matrices over Commutative $*$-Algebras: Quasi-Unitary Matrices}\label{sect_sos_com_matr}
In this section $\cA$ is a finitely generated unital commutative $*$-algebra without zero divisors over the field $\dL,$ where $\dL$ is $\dR$ or $\dC$. Let $\widehat{\cA}$ denote the set of characters on $\cA$, that is, $\widehat{\cA}$ is the set of all nontrivial $*$-homomorphisms $\chi:\cA\to\dL$. We assume that $\widehat{\cA}$ separates the elements of $\cA$. the latter implies in particular that $\sum_{j=1}^n a_j^*a_j^{}=0$ always implies that $a_1=\dots=a_n=0$. Define 
$$
\cA_+=\{ a\in \cA: \chi(a)\geq 0\ {\rm for}\ \chi \in \widehat{\cA}~\}~,~M_n(\cA)_+=\{ A{=}(a_{ij})\in M_n(\cA): (\chi(a_{ij}))\geq 0\ {\rm for}\ \chi \in \widehat{\cA}~\}.
$$
% without zero divisors for which Artins theorem holds (genauer aufschreiben!). 
Artin's theorem for the matrix algebra $M_n(\dR[t_1,\dots,t_d])$ was proved independently in \cite{gr} and in \cite{ps}. A constructive proof based on Schur complements was first developed in \cite{s4}. The aim of this short section is to give a proof of Artin's theorem for the matrix $*$-algebra $M_n(\cA)$ by using quasi-unitary matrices.

A matrix $T\in M_n(\cA)$ is called {\it quasi-unitary} if there exist an $s\in \sum\cA^2$ such that $T^*T= TT^*= sI$. Obviously, the element $s$ is uniquely determined by $T$. It will be denoted by $s(T)$. 

\begin{exam}
If $a_1,a_2\in \cA$, then
\begin{gather}\label{ext}
T=\left(
\begin{array}{ll}
a_1 & -a_2^* \\ 
a_2 & ~~a_1^* 
\end{array}\right)
\end{gather}
is a quasi-unitary matrix and $s(T)=a_1^*a_1+a_2^*a_2$.

To explain the approach in a simple special case we consider a $2\times 2$-matrix
\begin{gather*}
A=\left(
\begin{array}{ll}
a & b \\ 
b^* & c 
\end{array}\right)
\end{gather*}
of $M_2(\cA)_+$ such that $a\not= 0$. Then $a \in \cA_+$ and $ac-b^*b \in \cA_+$. Let us assume that there exist elements $a_1,a_2, d \in \cA$ such that $a=a_1^*a_1+a_2^*a_2$ and $ac-b^*b=d^*d$. 

Setting $(y_1,y_2)^t=T(b,d)^t$, that is, $y_1=a_1b-a_2^*d, y_2=a_2b+a_1^*d$, we have
\begin{gather*}
a^2A= \left(
\begin{array}{ll}
aa_1 & y_1 \\ 
aa_2 & y_2 
\end{array}\right)^*\left(
\begin{array}{ll}
aa_1 & y_1 \\ 
aa_2 & y_2 
\end{array}\right).
\end{gather*}
\end{exam}

\begin{prop}\label{quasit}
Let $k\in \dN$ and let $a_1,\dots,a_{2^k} \in \cA$ be given. There exists a quasi-unitary matrix $T{=}(t_{ij})\in M_{2^k}(\cA)$ and a nonzero element $s\in \sum \cA^2$ such that $t_{i1}=sa_i$ for $i=1,\dots,2^k$.
\end{prop}
\begin{proof} 
We proceed by induction on $k$. In the case $k=1$ we take the matrix $T$ given by (\ref{ext}). Suppose that the assertion is valid for $k$. If all elements $a_{2^k+1},\dots,a_{2^{k+1}}$ are zero, the assertion holds by the induction hypthesis. Suppose now that not all of these elements are zero. Then there are quasi-unitary matrices $A=(a_{ij}),B=(b_{ij}) \in M_{2^k}(\cA)$ and nonzero elements $s_1,s_2 \in \sum \cA^2$ such that $a_{i1}=s_1a_i$ and $b_{ij}=s_2 a_{2^k+i}$ for $i=1,\dots,2^k$. Since not all elements $a_{2^k+1},\dots,a_{2^{k+1}}$ are zero, we have $s(B)=s_2^2 \sum_{i=1}^{2^k} a_{2^k+i}^* a_{2^k+i}\not= 0$. A straightforward computation shows that the block matrix
 \begin{gather}\label{mata}
T=\left(
\begin{array}{ll}
s_2s(B)A & -s_1s(B)B^* \\ 
s_1s(B) B & ~~s_2BA^*B^* 
\end{array}\right)
\end{gather}
is quasi-unitary and we have $s:=s_1s_2s(B)\in \sum \cA^2$, $s\not= 0$ and $t_{i1}=sa_i$ for $i=1,\dots,2^{k+1}$.
%, where $s:=s_1s_2s(B)$.
\end{proof}

%The preceding proposition can be used to give another proof of Artins theorem for matrices %over commutative $*$-algebras. 
Artin's theorem for the $*$-algebra $M_n(\cA)$ is the following result.
\begin{thm}\label{artinmat}
Retain the assumptions and the notation stated above.
For any $A\in M_n(\cA)_+$ there exists an element $c \in \cA^\circ$ such that $c^*cA \in \sum M_n(\cA)^2$.
\end{thm}
\begin{proof}
Assume that the assertion holds for $n \in \dN$. Let $A\in M_{n+1}(\cA)_+$. We write $A$ as a block matrix
\begin{gather}\label{mata}
A=\left(
\begin{array}{ll}
a & \gb \\ 
\gb^* & C
\end{array}\right),
\end{gather}
where $a\in\cA,\gb\in M_{1,n}(\cA)$,
%\ \gb^*:=(b_1^*,\dots,b_n^*)^t\in M_{n,1}(\cA)$ 
and $C\in M_n(\cA)$. If $a=0$, the determinants of all $2 \times 2$ principal submatrices containing $a$ are in $\cA_+$ (see e.g. \cite{zh}, p. 161) which in turn implies that $\gb=0$. (This can be also derived from Corollary \ref{cor_1}.) Then the assertion follows a by applying the induction hypothesis to $C$. From now on we suppose that $a\not= 0$.

Since $A\in M_n(\cA)_+$, we have $a \in \cA_+$ and $D:=aC-\gb^*\gb \in M_{n}(\cA)_+$ (see e.g. \cite{zh}, p. 184). We assumed that $\cA$ is a finitely generated unital commutative$*$-algebra such that $\widehat{\cA}$ separates the elements of $\cA$. Therefore, Artin's theorem holds for $\cA$ which also gives the first induction step. By this theorem and by the induction hypothesis there exist elements $c_1,c_2 \in \cA^\circ$ such that 
$c_1^*c_1a \in \sum \cA^2$ and $c_2^*c_2D \in \sum M_n(\cA)^2$. 
Setting $c=c_1c_2$ and adding zeros in the sums of squares if needed, there are $m \in \dN$, elements $a_1,\dots,a_{2^m} \in \cA$ and column matrices $\ga_2,\dots,\ga_{2^m}\in M_{1,n}(\cA)$ such that 
$c^*ca =\sum_{j=1}^{2^m} a_j^*a_j$ and $c^*cD = \sum_{j=2}^{2^m} \ga_j^* \ga_j$. Setting $\ga_1:=c\gb$, we obtain
$c^*ca C=\sum_{j=1}^{2^m} \ga_j^* \ga_j$.

By Proposition \ref{quasit} there exist an $r\in \sum\cA^2$ and a quasi-unitary matrix $T=(t_{ij}) \in M_{2^m}(\cA)$ such that $t_{i1}=r a_i$ for $i=1,\dots,2^m$. Since $T$ is quasi-unitary, we then have $$s(T)=r^2 \sum_{j=1}^{2^m} a_j^*a_j =r^2c^*ca. $$
Since $r$, $c$ and $a$ are nonzero elements and $\cA$ has no zero divisors, $s(T)\not= 0$. Putting $\gy_j=\sum_{k=1}^m rt_{jk}\ga_k$ for $j=1,\dots,2^m$, we compute 
\begin{align*}
&\sum_{j} \gy_j^*\gy_j = \sum_j \sum_{k,l} r^2\ga_k t_{jk}^*t_{jl} \ga_l = 
 \sum_{k,l} r^2 \ga_k (T^*T)_{kl} \ga_l = r^2 s(T) \sum_k \ga_k^* \ga_k =r^2s(T)c^*c aC =s(T)^2 C,\\
&\sum_j (rcaa_j)^* \gy_j =\sum_{j,k} r^2c^*at_{j1}^* t_{jk} \ga_k =
r^2c^*a \sum_k (T^*T)_{1,k} \ga_k= r^2c^*a s(T) \ga_1 = r^2 c^*ca s(T)\gb =s(T)^2 \gb,\\
& \sum_{j} \gy_j^* rcaa_j =(\sum_j (rcaa_j)^*\gy_j)^* =s(T)^2 \gb^*,~~\sum_{j} (rcaa_j)^*rcaa_j =r^2c^*ca^2 \sum_j a_j^*a_j =s(T)^2 a.
\end{align*}
Setting $\gx_j:=(rcaa_j,\gy_j)\in M_{1,n+1}(\cA)$, the preceding four equations mean that
\begin{align*}
\sum_{j=1}^{2^m} \gx_j^*\gx_j = \sum_{j=1}^{2^m}
\left(
\begin{array}{ll}
(rcaa_j)^* & 0 \\ 
~~~\gy_j^* & 0
\end{array}\right) 
\left(
\begin{array}{ll}
rcaa_j & \gy_j \\ 
~~0 & 0
\end{array}\right)=
s(T)^2 \left(
\begin{array}{ll}
a & \gb \\ 
\gb^* & C
\end{array}\right) =s(T)^2 A.
%\end{gather}
\end{align*}
\end{proof}

\begin{remk}
The notion of quasi-unitaries used above might be useful for {\it non-commutative} $*$-algebras as well.
Let $\cB$ be a (not necessarily commutative) unital $*$-algebra with center $\cZ(\cB)$. Suppose that $bz=0$ for $b\in \cB$ and $z\in \cZ(\cB)$ always implies that $b=0$ or $z=0$. An element $b \in \cB$ is called a {\it quasi-unitary} if there exists an element $s(b)\in \sum \cZ(\cB)^2$
 such that $b^*b=bb^* =s(b).$
Clearly, if $a$ and $b$ are quasi-unitary, so is $ab$ and $s(ab)=s(a)s(b)$.
\end{remk}
\section{Matrix polynomials on intervals}\label{sect_sos_interval}

In this section we give another application of conditional expectations. Let $\dC[\dT]$ be the algebra of trigonometric polynomials in one variable with complex coefficients. The classical Fej\'er-Riesz theorem says that every positive element $f\in\dC[\dT]$ is equal to $g\overline{g}$ for some $g\in\dC[\dT].$ The following non-commutative version of this theorem follows from Theorem 7 in \cite{ro}.
\begin{prop}\label{prop_fej_rsz}
A self-adjoint element $X\in M_n(\dC[\dT])$ is positive semidefinite in every point $t\in\dT$ if and only if $X=Y^*Y$ for some $Y\in M_n(\dC[\dT]).$
\end{prop}

Let us identify the $*$-algebra $\dC[\dT]$ with the quotient $*$-algebra $\dC[x,y]/\langle 1-x^2-y^2\rangle.$ Since each element of $\dC[x,y]/\langle 1-x^2-y^2\rangle$ can be written uniquely in the form $f_1+yf_2$ with $ f_1,f_2\in\dC[x],$ under this identification $\dC[\dT]$ becomes the vector space $\dC[x]+y\dC[x]$ with multiplication rule $$(f_1+yf_2)(g_1+yg_2)=f_1f_2+(1-x^2)f_2g_2+y(f_1g_2+f_2g_1),\ f_i,g_i\in\dC[x].$$ We consider $\dC[x]$ as a $*$-subalgebra of $\dC[\dT]$ and define a mapping $\gp:\dC[\dT]\to\dC[x]$ by $\gp(f_1(x)+yf_2(x)):=f_1(x).$ In the same manner we consider $M_n(\dC[x])$ as a $*$-subalgebra of $M_n(\dC[\dT])$ and define a projection $\gP:M_n(\dC[\dT])\to M_n(\dC[x])$ by $\gP((f_{ij})):=(\gp(f_{ij})),\ f_{ij}\in\dC[\dT].$ 

\begin{prop}\label{prop_cexp_circ}
The mapping $\gP$ is a conditional expectation such that $$\gP(X^*X)=Y_1^*Y_1^{}+(1-x^2)Y_2^*Y_2^{},\ X\in\dC[\dT],\ Y_1,Y_2\in M_n(\dC[\dT]).$$
\end{prop}
\begin{proof}
We write $X$ as $X_1+yX_2,\ X_1,X_2\in M_n(\dC[x]).$ It is easily seen that $\gP(X^*X)=X_1^*X_1^{}+y^2X_2^*X_2^{}=X_1^*X_1^{}+(1-x^2)X_2^*X_2^{}.$ One can readily check that $\gP$ fulfils the other axioms of a conditional expectation.
\end{proof}

\begin{thm}\label{thm_matr_pol_int}
Let $F(x)\in M_n(\dC[x])$ be a self-adjoint matrix polynomial. Then $F(x)$ is positive semidefinite in every point $x\in[a,b],a<b$ if and only if $F=G_1^*G_1^{}+(b-x)(x-a)G_2^*G_2^{}$ for some $G_1,G_2\in M_n(\dC[x]).$
\end{thm}
\begin{proof}
It suffices to prove the only if direction. Upon applying a linear transformation we can assume that $[a,b]=[-1,1].$ Let $F(x)$ be positiveon $[-1,1].$ Then $F(x)$ is a positive element of $M_n(\dC[\dT])$, so that $F=G^*G$ for some $G\in M_n(\dC[\dT])$ by Proposition \ref{prop_fej_rsz}. From Proposition \ref{prop_cexp_circ} we obtain $F=\gP(F)=\gP(G^*G)=G_1^*G_1^{}+(1-x)(x+1)G_2^*G_2^{}$ for some $G_1,G_2\in M_n(\dC[x]).$
\end{proof}

The same methods can be used to derive a Positivstellensatz for matrix polynomials on $[a,+\infty).$ First we recall a well-known result (see e.g. \cite{dj}).
\begin{prop}\label{prop_matr_pol_1dim}
A self-adjoint matrix polynomial $X\in M_n(\dC[x])$ is positive semidefinite for all $x\in\dR$ if and only if $X=Y^*Y$ for some $Y\in M_n(\dC[x]).$
\end{prop}

Proceeding in the same way as above, we write $g\in\dC[x]$ uniquely as $g_1+xg_2,\ g_1,g_2\in\dC[x^2]$ and define a conditional expectation $\gp_2:\dC[x]\to\dC[x^2]$ by $\gp_2(g):=g_1$. We then obtain 
\begin{thm}\label{thm_matr_pol_hint}
Let $F(x)\in M_n(\dC[x])$ be a self-adjoint matrix polynomial. Then $F(x)$ is positive semidefinite for all $x\in[a,+\infty)$ if and only if $F=G_1^*G_1^{}+(x-a)G_2^*G_2^{}$ for some $G_1,G_2\in M_n(\dC[x]).$
\end{thm}

\noindent{\bf Remark.} For $F(x)=(F_{ij}(x))\in M_n(\dC[x])$ we define  $\deg F:=\max_{ij}\deg F_{ij}.$ A closer look at the proofs of  Propositions \ref{prop_fej_rsz} and \ref{prop_matr_pol_1dim} allows us to estimate the degrees of $G_1$ and $G_2$ in Theorems \ref{thm_matr_pol_int} and \ref{thm_matr_pol_hint}. In Theorem \ref{thm_matr_pol_int} we can achieve $\deg G_1\leq\deg F,\ \deg G_2\leq\deg F-1$ and in Theorem \ref{thm_matr_pol_hint} we can have $\deg G_1\leq \frac{1}{2}\deg F,\deg G_2\leq \frac{1}{2}(\deg F-1),$ where we mean $G_2=0$ when $\deg G_2=-1$.

%In the last section we discuss the problem of generalizing the above Theorems to other %semi-algebraic sets.
\section{Positivstellens\"atze for matrices over fields}\label{sect_sos_matr_field}

In this section $\bR$ is a formally real field and recall that $\bR_+$ denotes the set of finite sums of squares in $\bR.$ Let $\bC=\bR(\sqrt{-1})=\bR+\img\bR,$ where $\img=\sqrt{-1}$, with involution defined by $(l_1+\img l_2)^*=l_1-\img l_2,\ l_1,l_2\in\bR.$ The purpose of this section is to study positivity in $M_n(\bR)$ and $M_n(\bC)$ and to prove the Positivstellensatz stated as Theorem \ref{posfield} below. 

In what follows the symbol $\bL$ denotes one of the fields $\bR$ or $\bC$ meaning that a statement or definition holds for both $\bR$ and $\bC.$ We consider $M_n(\bL)$ as $*$-algebra with  involution $(a_{ij})^*=(a_{ji}^*)$.  We denote by $\gp_{kk}:M_n(\bL)\to\bL$ the mapping $(a_{ij})\mapsto a_{kk}$ and by $\tr:M_n(\bL)\to\bL$ the mapping $\tr=\frac{1}n(\gp_{11}+\dots+\gp_{nn}).$ 

Let $B=\diag{1,\lambda_1,\dots,\lambda_{n-1}}\in M_n(\bR)$ be a diagonal matrix, where $\lambda_1,\dots,\lambda_{n-1}\in\bR^\circ$, and  $\tau$ the involution on $M_n(\bL)$ associated with $B,$ that is, $X^\tau:=B^{-1}X^*B.$ Let $\langle\cdot,\cdot\rangle$ be the standard inner product on $\bL^n$  and  $\langle x,y\rangle_1:=\langle Bx,y\rangle$ the inner product on $\bL^n$ defined by $B.$ Then $\tau$ is the involution associated to $\langle\cdot,\cdot\rangle_1.$ We first record a simple lemma.
%Our main result is a Positivstellensatz for $(M_n(\bL),\tau).$ 

\begin{lemma}\label{lemma_cond_star_ord}
Let $p$ be an ordering on $\bR.$ Then the following are equivalent:
\begin{enumerate}
 \item[\rm (i)] for  fixed $k\in\set{1,\dots,n}$ and for all $X\in M_n(\bL),$ $p$ contains $\gp_{kk}(X^\tau X),$
 \item[\rm (ii)] for all $k\in\set{1,\dots,n}$ and for all $X\in M_n(\bL),$ $p$ contains $\gp_{kk}(X^\tau X),$
 \item[\rm (iii)] for all $X\in M_n(\bL),$ $p$ contains $\tr(X^\tau X),$
 \item[\rm (iv)] for all $\gf\in\bL^n,$ $p$ contains $\langle\gf,\gf\rangle_1.$
 \item[\rm (v)] $p$ contains $\lambda_1,\dots,\lambda_{n-1}.$
\end{enumerate}
\end{lemma}
\begin{proof}
Equivalence (i)$\Leftrightarrow$(ii)$\Leftrightarrow$(iii)$\Leftrightarrow$(iv) follows from simple computations. Let $\gee_1,\dots\gee_n$ be the standard base of $\bL^n$ and let $\gf=\sum_{i=1}^n f_i\gee_i.$ Then $$\langle\gf,\gf\rangle_1=\sum_{i=1}^n f_i^*f_i^{}\langle\gee_i,\gee_i\rangle_1=\sum_{i=1}^n \lambda_{i-1}^{}f_i^*f_i^{},\ \mbox{where}\ \lambda_0=1$$ which implies (iv)$\Leftrightarrow$(v).
\end{proof}

\begin{defn}\label{defn_star_ord_Mn}
An ordering $p$ of $\bR$ is called $*$-\textit{ordering} if one of the statements (i)-(v) in Lemma \ref{lemma_cond_star_ord} is satisfied.
\end{defn}

To avoid degenerate situations, we assume that $*$-orderings exist. We denote by $P_B(\bL)$ the preordering generated by $\lambda_1,\dots,\lambda_{n-1}.$ The existence of $*$-orderings is equivalent to the requirement that $P_B(\bL)$ is proper, that is, $-1\notin P_B(\bL)$. Since a proper preordering is the intersection of all orderings containing it, $P_B(\bL)$ is the intersection of all $*$-orderings. An element $l\in\bL$ is positive in all $*$-orderings if and only if $l$ belongs to $P_B(\bL).$

% The next proposition is a standard statement from linear algebra in the case $\bL=\dR$ and $B=I,$ see e.g. Theorem 7.4., \cite{ip} an its proof is completely analogous.
% \begin{prop}\label{prop_jac_meth}
% Let $X=X^\tau\in M_n(\bL).$ Then there exists $Y=(y_{ij})\in M_n(\bL)$ satisfying $y_{ii}=1,\ y_{ij}=0,\ j>i$ such that $YXY^\tau$ is diagonal. If the leading principal minors $\Delta_1,\dots,\Delta_n$ of $X$ are non-zero then 
% \begin{gather}
% YXY^\tau=\diag{\Delta_1,\ \frac{\Delta_2}{\Delta_1},\dots, \frac{\Delta_n}{\Delta_{n-1}}}
% \end{gather}
%\end{prop}

\begin{lemma}\label{lemma_pos_in_Mn_K}
For $X=X^\tau\in M_n(\bL)$ the following are equivalent:
\begin{enumerate}
 \item[\rm (i)] for a fixed $k\in\set{1,\dots,n}$ and every $Y\in M_n(\bL)$, we have $\gp_{kk}(Y^\tau XY)\in P_B(\bL),$
 \item[\rm (ii)] for every $k\in\set{1,\dots,n}$ and every $Y\in M_n(\bL)$, we have $\gp_{kk}(Y^\tau XY)\in P_B(\bL),$
 \item[\rm (iii)] for every $Y\in M_n(\bL)$, we have $\tr(Y^\tau XY)\in P_B(\bL),$
 \item[\rm (iv)] for every $\gf=(f_1,\dots,f_n)\in\bL^n$, we have $\langle X\gf,\gf\rangle_1\in P_B(\bL),$
 \item[\rm (v)] all of the principal minors of $X$ belong to $P_B(\bL).$
\end{enumerate}
\end{lemma}
\begin{proof}
The equivalence of the statements (i)-(iv) is proved by simple computations. The equivalence of (iv) and (v) is a slight generalization of Sylvester's criterion for positive semi-definite matrices.
\end{proof}

\begin{defn}\label{defn_star_ord_Mn}
Let $(M_n(\bL),\tau)_+$ be the set of all $X{=}X^\tau\in M_n(\bL)$ for which one of the statements (i)-(v) in Lemma \ref{lemma_pos_in_Mn_K} is satisfied.
Such elements are called positive (with respect to $\tau$).
\end{defn}

%Denote by $(M_n(\bL),\tau)_+$ the set of all positive elements in $M_n(\bL).$ 
The following Lemma contains some elementary properties of $(M_n(\bL),\tau)_+.$

\begin{lemma}
%The following statements hold:
\begin{enumerate}
	\item[\rm (i)] $\sum X_i^\tau X_i^{}\in (M_n(\bL),\tau)_+$ for arbitrary $X_i\in M_n(\bL).$
	\item[\rm (ii)] If $Y=Y^\tau\in(M_n(\bL),\tau)_+$ then $X^\tau YX^{}\in(M_n(\bL),\tau)_+$ for every $X\in M_n(\bL).$
\end{enumerate}
\end{lemma}

The following proposition seems to be missing in the literature even for the case when $B$ is an identity matrix and $\bL=\bR.$
\begin{prop}\label{prop_prod_com_pos} 
If $X,Y\in (M_n(\bL),\tau)_+$  and $XY=YX,$ then $XY\in (M_n(\bL),\tau)_+.$ In particular, $(M_n(\bL),\tau)_+$ is a non-commutative preordering.
\end{prop}
\begin{proof}
We give the proof for the case $\bL=\bR.$ The case $\bL=\bC$ is treated similarly. Let us fix a $*$-ordering $p$ on $\bR$ and let $\overline{\bR}$ be the real closure of $(\bR,p).$ It follows from Lemma \ref{lemma_pos_in_Mn_K} (v) that $X,Y\in (M_n(\overline{\bR}),\tau)_+$ and it is enough to show that $XY\in (M_n(\overline{\bR}),\tau)_+.$ %This is done by a simultaneous digonalization of $X$ and $Y$ in $M_n(\overline{\bR}).$ 
Since $\overline{\bR}$ is real closed and $p$ contains $\lambda_1,\dots,\lambda_{n-1},$ the matrix $B^{1/2}:=\diag{1,\sqrt{\lambda_1},\dots,\sqrt{\lambda_{n-1}}}$ belongs to $M_n(\overline{\bR}).$ Then the mapping $\phi(X)=B^{1/2}XB^{-1/2}$ defines a $*$-isomorphism of $(M_n(\overline{\bR}),\tau)$ and $M_n(\overline{\bR})$ with respect to the transpose involution. The elements $X_1:=\phi(X)$ and $Y_1:=\phi(Y)$ belong to $M_n(\overline{\bR})_+$ and commute. It is enough to check that $X_1Y_1\in M_n(\overline{\bR})_+.$ This is proved by simultaneous diagonalization of $X_1$ and $Y_1.$ Since $X_1$ is symmetric, there exists an orthogonal matrix $U\in M_n(\overline{\bR})$ such that $U^TX_1U$ is a diagonal matrix $\diag{x_1,\dots,x_n}\in M_n(\overline{\bR}).$ Assume that the $x_i$ are pairwise different (the other case is treated by choosing $U$ more specifically). Then $U^TY_1U$ commutes with $U^TX_1U$ and since the $x_i$ are pairwise different, $U^TY_1U$ is also diagonal. Since $U$ is orthognoal, both $U^TX_1U$ and $U^TY_1U$ belong to $M_n(\overline{\bR})_+.$ Hence their product $U^TX_1Y_1U$ also belongs to $M_n(\overline{\bR})_+$ and $X_1Y_1\in M_n(\overline{\bR})_+.$  
\end{proof}

\begin{exam}\label{exam_ctr_ps}
It is proved in \cite{ps} that for $n=2$ we have $(M_n(\bR),\tau)_+=\sum(M_n(\bR),\tau)^2$ for arbitrary $\bR.$ For $n=3$ there is the following counterexample given in \cite{klep}. Let $\bR=\dR(s,t)$ be the field of rational functions in two variables and $B=\diag{s,t,st}.$ Then $\diag{st,st,st}$ is positive, but it is not a sum of squares in $(M_n(\bR),\tau),$ see Theorem 3.2 in \cite{klep}.
\end{exam}

Let $P_B^n(\bL)\subseteq M_n(\bR)$ denote the set of all diagonal matrices with entries from $P_B(\bL).$ 
\begin{lemma}\label{lemma_pos_diag}
Let $X\in M_n(\bR)$ be a diagonal matrix. Then $X\in(M_n(\bR),\tau)_+$ if and only if $X\in P_B^n(\bL).$
\end{lemma}
\begin{proof}
This follows immediately from Lemma \ref{lemma_pos_in_Mn_K}(iv).
\end{proof}

The next proposition is a standard fact from linear algebra in the case $\bL=\dR$ and $B=I.$ Its proof in the present case is completely analogous.
\begin{prop}\label{prop_jac_meth}
Let $X=X^\tau\in M_n(\bL).$ Then there exists an invertible matrix $Y=(y_{ij})\in M_n(\bL)$ such that $Y^\tau XY$ is diagonal. 
\end{prop}

By Example \ref{exam_ctr_ps}, $P_B^n(\bL)$ contains elements which are not sums of squares. We prove that $P_B^n(\bL)$ generates $(M_n(\bR),\tau)_+$ as a quadratic module. 
\begin{prop}\label{prop_pst_alg_inv}
Let $X=X^\tau\in M_n(\bL).$ Then $X\in(M_n(\bL),\tau)_+$ if and only if there exist $Y\in M_n(\bL)$ and $D\in P_B^n(\bL)$ such that $X=Y^\tau DY.$
\end{prop}
\begin{proof}
First let $X=Y^\tau DY$ and $ D\in P_B^n(\bL).$ Then $\langle Xv,v\rangle_1=\langle Y^\tau DYv,v\rangle_1=\langle DYv,Yv\rangle_1\geq 0$ in all $*$-orderings.

Conversely, suppose that $X\in(M_n(\bL),\tau)_+$. By Proposition \ref{prop_jac_meth} there exists an invertible matrix $Z\in M_n(\bR)$ such that $D:=ZXZ^\tau$ is diagonal. Since $X\in (M_n(\bR),\tau)_+$ we also have $ZXZ^\tau\in (M_n(\bR),\tau)_+.$ Since $D$ is diagonal, Lemma \ref{lemma_pos_diag} implies that $D\in P_B^n(\bL).$ Therefore, $X=Z^{-1}D(Z^\tau)^{-1}=Z^{-1}D(Z^{-1})^\tau.$ Put $Y=(Z^{-1})^\tau.$
\end{proof}

\begin{prop}\label{prop_P_B_ncsos}
Every element of $P_B^n(\bL)\subseteq(M_n(\bR),\tau)_+$ is a non-commutative sum of squares in $(M_n(\bR),\tau).$ In the case $n=2$ we have $P_B^n(\bL)\subseteq\sum(M_n(\bR)\tau)^2.$
\end{prop}
\begin{proof}
It is enough to prove the assertion for the matrices $\lambda_m\cdot E_{kk}\in P_B(\bL).$ For $k>1$ we compute
$$(\lambda_{k-1}E_{m+1,k})^\tau(\lambda_{k-1}E_{m+1,k})=\lambda_{k-1}^2 B^{-1}E_{k,m+1}BE_{m+1,k}=\lambda_m\lambda_{k-1} E_{kk},$$
and we also have $E_{1k}^\tau E_{1k}^{}=\lambda_{k-1}^{-1} E_{kk}.$ Thus, $\lambda_m\cdot E_{kk}$ is the product of two commuting squares $(\lambda_{k-1}E_{m+1,k})^\tau(\lambda_{k-1}E_{m+1,k})$ and $E_{1k}^\tau E_{1k}^{}.$ In the case $k=1$, $\lambda_m\cdot E_{kk}=(\lambda_m E_{m+1,1})^\tau(\lambda_m E_{m+1,1}).$

For $n=2$ the inclusion $P_B^n(\bL)\subseteq\sum(M_n(\bR)\tau)^2$ follows from $\lambda_1E_{11},\lambda_1E_{22}\in\sum(M_n(\bR)\tau)^2.$
\end{proof}

\begin{thm}\label{posfield}
Let $X=X^\tau\in M_n(\bL)$. Then $X\in (M_n(\bL),\tau)_+$ if and only if it is a non-commutative sum of squares. For $n=2$ we have the equality $(M_n(\bR),\tau)_+=\sum(M_n(\bR)\tau)^2.$
\end{thm}
\begin{proof} The assertion follows from Propositions \ref{prop_prod_com_pos}, \ref{prop_pst_alg_inv} and \ref{prop_P_B_ncsos}.
\end{proof}

\noindent\textbf{Remarks.} 1. The condition $-1\notin P_B(\bL)$ seems to be the natural definition for $(M_n(\bR),\tau)$ of being formally real. It implies that an equality $\sum_iX_i^\tau X_i^{}=0$ is only possible when all $X_i=0.$
\medskip

\noindent 2. The case of a non-trivial involution on the basic field can be treated more generally. Then one defines $\bC$ as $\bR(\sqrt{-s})$ for some $s\in\sum\bR^2,\ s\neq 0,$ that is, one takes $\mathbf{j}=\sqrt{-s}$ as imaginary unit. Since $\bR$ is formally real, $\bC$ is a proper extension of $\bR$ and it is not formally real. A natural involution on $\bC$ is defined by $(l_1+\mathbf{j}l_2)^*=l_1-\mathbf{j}l_2.$ A field $\bC$ arising in this way might be called "formally complex." The involution on $\bC$ has the property that sums of "hermitian squares" $l_i^*l_i^{}$ are contained in $\sum\bR^2.$ All definitions and statements of this section carry over to this case.

\section{Positivstellens\"atze for cyclic algebras}\label{sect_sos_cyc}
Let $\bR,\bC$ and $\bL$ be as in the preceding section. Suppose that $\bL/\bK$ is a Galois extension with group $\cyc{n}$, $\sigma$ is an automorphism of $\bL$ which generates $\cyc{n},$ and $\gA$ is a cyclic algebra associated to $\bL/\bK.$ Our aim in this section is to develop a Positivstellensatz (Theorem \ref{thm_pos_cyc_alg}) for the cyclic algebra $\gA$.

Since $\gA$ is a cyclic algebra, there exist fixed elements $e\in\gA$ and $\ a\in \bK^{\circ}$ such that
\begin{gather}\label{eq_gA_cyclic}
 \gA=\Un{}\cdot\bL\oplus e\cdot\bL\oplus\dots e^{n-1}\cdot\bL,\ e^n=a\cdot\Un{},\ \mbox{and}\ \lambda\cdot e=e\cdot\sigma(\lambda),\ \mbox{for}\ \lambda\in \bL.
\end{gather}

From (\ref{eq_gA_cyclic}) it follows that $\gA=(\bL/\bK,\sigma,a)$ is a $\cyc{n}$-graded algebra with standard grading $\gA_k=e^k\cdot\bL.$ We assume in addition that $\gA$ is a $\cyc{n}$-graded $*$-algebra, that is, the involution of $\gA$ satisfies 
\begin{gather}\label{eq_gA_grad_star}
\gA_k^*=\gA_{-k}^{}\ \mbox{for}\ k\in\cyc{n}. 
\end{gather}

By (\ref{eq_gA_grad_star}), $\gA_0\simeq\bL$ is invariant under the involution, where the involution of $\bL$ is defined as in the previous section. It follows from (\ref{eq_gA_grad_star}) that $e^*\in\gA_{n-1},$ so there exists an element $l_0\in\bL^{\circ}$ such that $e^*=l_0\cdot e^{n-1}.$ Hence $e^*e=l_0\cdot e^n=l_0\cdot a.$

\begin{lemma}\label{lemma_sigma*}
 $\sigma$ is a $*$-automorphism of $\bL$. 
\end{lemma}
\begin{proof}
Applying the involution to the equality $e\cdot\sigma(l)=l\cdot e$ and replacing $e^*$ by $l_0\cdot e^{n-1}$ we obtain:
$$
\sigma(l)^*l_0\cdot e^{n-1}=l_0\cdot e^{n-1}\cdot l^*=l_0\sigma(l^*)\cdot e^{n-1}=\sigma(l^*)l_0\cdot e^{n-1}.
$$
This implies that $\sigma(l)^*=\sigma(l^*).$
\end{proof}

%It follows from (A1) that $le^k=\sigma^k(l)e^k,$ which implies the following formula:
%\begin{gather}\label{eq_sigma_l}
%\sigma^k(l)=\frac{e^{k*}le^k}{e^{k*}e^k},\ \mbox{for all}\ l\in\bL.
%\end{gather}

%Proposition ?? below gives the connection between Definition \ref{defn_star_ord_Mn} and Definition \ref{defn_pos_in_gA}.

%\begin{remk} 
%It seems to us, that condition (A4) gives the most natural definition of $\gA$ to be {\it formally real}, cf. \cite{cim}. 
%\end{remk}

%Definition \ref{defn_pos_in_gA} can be expressed in terms of a regular representation of $\gA.$ 

Let $\epsilon$ be the left regular representation of $\gA.$ 
%More precisely, consider the base $\gee_k=e^{k-1},\ k=\overline{1,n}$ of $\gA$ as of the right $\bL$-vector space. 
We identify $\gA$ with $\bL^n$ via (\ref{eq_gA_cyclic}) by considering $\gee_k:=e^{k-1}$, $\ k=1,\dots,n$, as standard base in $\bL^n.$ Then the action $\epsilon$ of $\gA$ on $\bL^n$ is given by
\begin{gather}\label{eq_reg_rep_gA}
\epsilon(l)\gee_k=le^{k-1}=e^{k-1}\sigma^{k-1}(l)=\gee_k\sigma^{k-1}(l),\ \mbox{and}\\ 
\nonumber\epsilon(e)\gee_k=ee^{k-1}=e^k=\left\{\begin{array}{ll}\gee_{k+1}, & \mbox{if}\ k<n; \\ a\gee_1, & \mbox{for}\ k=n,\end{array}\right\}
\end{gather}
That is, these formulas define an algebra homomorphism $\epsilon:\gA\to M_n(\bL)$. From (\ref{eq_gA_cyclic}) it follows that $\epsilon$ is injective and that the corresponding matrices of $\epsilon(l)$ and $\epsilon(e)$ are
\begin{gather}\label{eq_eps_matr_form}
\epsilon(l)=\left(
\begin{array}{ccccc}
    l   	&  0		&   \dots	&   		&  0  			\\
    0   	& \sigma(l) 	&   		&   		&  \vdots   		\\
    \vdots    	& 	 	& \sigma^2(l)  	&   		&     			\\
        	&   		&  		& \ddots  	&     			\\
        	&   		&   		&   		&  0  			\\
    0		& \dots		&   		&  0 		&    \sigma^{n-1}(l)
\end{array}
\right),\ \epsilon(e)=\left(
\begin{array}{ccccc}
  0 		& \dots		&   			& \dots	 & a      	\\
  1 		& 0 			&   			&   		 &      						\\
    		& 1				& \ddots 	&   		 &      						\\
    		&   			& \ddots 	&     	 &      						\\
    		&   			&   			& 			 & 							\\
    		&   			&   			& 		 1 & 0
\end{array}
\right)
\end{gather}

We denote by $\gp$ the canonical projection from $\gA$ onto $\bL,$ that is, $\gp(\sum e^kl_k):=l_0.$ Note that $\gp$ maps self-adjoint elements of $\gA$ into $\bR.$ From (\ref{eq_eps_matr_form}) we easily derive the following 
\begin{lemma}\label{lemma_gp_gp11}
$\gp$ coincides with $\gp_{11}(\epsilon(\cdot)),$ where $\gp_{11}:M_n(\bL)\to\bL$ is defined in Section \ref{sect_sos_matr_field}.
\end{lemma}

Let $\langle\cdot,\cdot\rangle_1:\gA\times\gA\to\bL$ be a $\bL$-valued "inner product" on $\gA$ defined by
\begin{gather}\label{eq_gp_y*x}
\langle x,y\rangle_1:=\gp(y^*x),\ x,y\in\gA.
\end{gather}

For $x=\sum_k e^kx_k$ and $\ y=\sum e^my_m$, using (\ref{eq_gA_grad_star}) we compute
\begin{gather}\label{eq_sc_prod_gp}
\langle x,y\rangle_1=\sum_{k,m}\gp(y_m^*e^{m*}e^kx_k^{})=\sum_k (e^{*k}e^k)y_k^*x_k^{}.
\end{gather}

Set $B=\diag{\lambda_0,\lambda_1,\dots,\lambda_{n-1}},$ where $\lambda_k:=e^{*k}e^k,\ k=1,\dots,n-1,\ \lambda_0=1$ and let $P_B(\bL),\ P_B^n(\bL)$ be as in the preceding section. As above we assume that the space of $*$-orderings of $\bR$ is not empty, i.e. $-1\notin P_B(\bL).$ 

\begin{defn}\label{defn_pos_in_gA}
$\gA_+$ is the set of elements $y=y^*\in\gA$ such that $\gp(x^*yx)\in P_B(\bL)$ for all $x\in \gA$. The elements of $\gA_+$ are called positive. 
%\epsilon(e)
\end{defn}
Let $\langle x,y\rangle_1=\langle Bx,y\rangle$ and the involution $\tau$ on $M_n(\bL)$ be as in Section \ref{sect_sos_matr_field}.

\begin{prop}\label{prop_eps_*hom}
The embedding $\epsilon:\gA\to(M_n(\bL),\tau)$ is a $*$-homomorphism. Moreover, $\epsilon(\gA_+)=\epsilon(\gA)\cap(M_n(\bL),\tau)_+.$
\end{prop}
\begin{proof}
The first statement follows from the computation
$$
\langle\epsilon(z)x,y\rangle_1=\gp(y^*zx)=\gp((z^*y)^*x)=\langle x,\epsilon(z^*)y\rangle_1.
$$
Let $y=y^*\in\gA_+$ and $x\in\gA\simeq\bL^n.$ Then $\gp(x^*yx)\in P_B(\bL)$ if and only if $\langle\epsilon(y)x,x\rangle_1\in P_B(\bL).$ Thus $\epsilon(y)\in(M_n(\bL),\tau)_+$ if and only if $y\in\gA_+.$
\end{proof}

We identify $\gA$ with the $*$-subalgebra $\epsilon(\gA)\subseteq M_n(\bL).$
	
\begin{lemma}\label{lemma_PB_sigma_inv}
$P_B(\bL)\subseteq\bL$ is invariant under $\sigma.$
\end{lemma}
\begin{proof}It follows from (\ref{eq_gA_cyclic}) that $(e^{*}e)\sigma(l)={e^{*}le}\ \mbox{for all}\ l\in\bL.$ Using this fact we get
$$\sigma(e^{*k}e^k)=e^{*(k+1)}e^{k+1}(e^*e)^{-1}=e^{*(k+1)}e^{k+1}\cdot e^*e\cdot((e^*e)^{-1})^2\in P_B(\bL).$$
\end{proof}

%Denote by $QM$ the quadratic module generated by $P_B(\bL).$ That is $\cC$ consists of the elements:
%\begin{gather}\label{eq_quadr_mod}
%\sum_i x_i^*y_ix_i,\ \mbox{where}\ x_i\in\gA, 
%\end{gather}
%and each $y_i$ is a product of a finite number of the elements (\ref{eq_set_e*e}). Our aim is to prove, that $\gA_+$ coincides with $\cC.$

It is known from the general theory that $\bL$ is a splitting field of $\gA$ and $\gA\otimes_\bK\bL\simeq M_n(\bL).$ Further, the Galois group $Gal(\bL/\bK)$ acts naturally by automorphisms on $M_n(\bL)$ such that $\gA$ is equal to the subalgebra of stable elements. Thus, the average of this action is a projection from $M_n(\bL)$ onto $\gA.$ 
%A projection $\gP:M_n(\bL)\to\gA$ which we construct below is of that kind. 
For $E_{mk}\otimes l\in M_n(\bL),\ l\in\bL$ we define 
\begin{gather*}
\gP(E_{mk}\otimes l):=\frac{1}n e^{m-k}\sigma^{-k+1}(l)=\frac{1}n e^{m-1}le^{-k+1}.
\end{gather*}

% Direct computations prove the following proposition (see also Proposition \ref{prop_homog_main} in Appendix).
\begin{prop}\label{prop_gP_bim_proj}
$\gP$ is a well-defined $\gA$-bimodule projection from $(M_n(\bL),\tau)$ onto $\gA.$ 
\end{prop}
\begin{proof}
We prove e.g. that $\gP$ preserves the involution. The other conditions follows by a similar reasoning. For $E_{mk}\otimes l\in M_n(\bL)$ we compute
\begin{gather*}
\gP((E_{mk}\otimes l)^\tau)=\gP(B^{-1}(E_{km}\otimes l^*)B)=\gP(E_{km}\otimes \lambda_{k-1}^{-1}\lambda_{m-1}^{}l^*)=\frac{1}n e^{k-1}\lambda_{k-1}^{-1}l^*\lambda_{m-1}^{}e^{-m+1}=\\
=\frac{1}ne^{k-1}\left(e^{*(k-1)}e^{k-1}\right)^{-1}l^*e^{*(m-1)}e^{m-1}e^{-m+1}=\frac{1}ne^{*(-k+1)}l^*e^{*(m-1)}=\\
=\frac{1}n(e^{m-1}le^{-k+1})^*=(\gP(E_{mk}\otimes l))^*.
\end{gather*}
\end{proof}

Let $QM_{M_n(\bL),\tau}(P_B^n(\bL))$ denote the quadratic module in $(M_n(\bL),\tau)$ generated by $P_B^n(\bL).$ By Proposition \ref{prop_pst_alg_inv} we have $\gA_+=QM_{M_n(\bL),\tau}(P_B^n(\bL)).$ Further, let $QM_\gA(P_B(\bL))$ be the quadratic module in $\gA$ generated by $P_B(\bL)\subseteq\bL\simeq\gA_0.$ The next proposition is a crucial step in the proof of the Positivstellensatz below.
\begin{prop}\label{prop_gP_cyc_alg_main}
$\gP$ maps $QM_{M_n(\bL),\tau}(P_B^n(\bL))$ onto $QM_\gA(P_B(\bL)).$
\end{prop}
\begin{proof}
Let $D=\diag{d_1,\dots,d_n}\in P_B^n(\bL),\ d_i\in P_B(\bL)$ and let $X=\sum_{ij}E_{ij}\otimes l_{ij}$ be an element from $M_n(\bL).$ The following computation is very similar to that used in the proof of Lemma \ref{lemma_sos_rank1}:
\begin{gather}
\nonumber X^\tau DX=\left(\sum_{i,j}E_{ij}\otimes l_{ij}\right)^\tau D\left(\sum_{k,m}E_{km}\otimes l_{km}\right)=B^{-1}\left(\sum_{i,j}E_{ji}\otimes l_{ij}^*\right)BD\left(\sum_{k,m}E_{km}\otimes l_{km}\right)=\\
\label{eq_sos_inv_r1}=\left(\sum_{i,j}E_{ji}\otimes l_{ij}^*\lambda_{j-1}^{-1}\lambda_{i-1}^{}d_i^{}\right)\left(\sum_{k,m}E_{km}\otimes l_{km}\right)=\sum_k\left(\sum_{j,m}E_{jm}\otimes l_{kj}^*l_{km}^{}\lambda_{j-1}^{-1}\lambda_{k-1}^{}d_k^{}\right)=\\
\nonumber=\sum_k\left(\sum_{j}E_{kj}\otimes l_{kj}\right)^\tau (E_{kk}\otimes d_k)\left(\sum_{m}E_{km}\otimes l_{km}\right)=\sum_kX_k^\tau DX_k^{},
\end{gather}
where $X_k=\sum_{m}E_{km}\otimes l_{km}.$ We now prove that for every $k=1,\dots,n$
\begin{gather}\label{eq_main_cyc}
\gP(X_k^\tau DX_k)=n^2\cdot\gP(X_k^\tau)\gP(E_{kk}\otimes d_k)\gP(X_k^{})\equiv n\cdot\gP(X_k)^*\sigma^{1-k}(d_k)\gP(X_k).
\end{gather}

Fix $k$. Using (\ref{eq_sos_inv_r1}) we compute
\begin{gather*}
\gP(X_k^\tau DX_k^{})=\gP\left(\sum_{j,m}E_{jm}\otimes l_{kj}^*l_{km}^{}\lambda_{j-1}^{-1}\lambda_{k-1}^{}d_k^{}\right)=\frac{1}n\sum_{j,m}e^{j-m}\sigma^{-m+1}(l_{kj}^*\lambda_{j-1}^{-1}\lambda_{k-1}^{}l_{km}^{}d_k^{}).
\end{gather*}
On the other hand, using the equality $(\gP(Y))^*=\gP(Y^\tau)=\gP(B^{-1}Y^*B)$ we derive
%
%By Proposition \ref{prop_gP_bim_proj} we have for all $Y\in M_n(\bL).$ Using this, we do the following computation in $\gA$ for a fixed $k=1,\dots,n:$
\begin{gather*}
\left(\gP(X_k)\right)^*\gP(E_{kk}\otimes d_k)\gP(X_k)%=\frac{1}{n}\gP\left(\left(\sum_{j}E_{kj}\otimes l_{kj}\right)^\tau\right)\cdot\sigma^{-k+1}(d_k)\cdot\gP\left(\sum_{m}E_{km}\otimes l_{km}\right)=\\
=\frac{1}{n}\gP\left(\sum_{j}E_{jk}\otimes l_{kj}^*\lambda_{j-1}^{-1}\lambda_{k-1}^{}\right)\cdot\sigma^{-k+1}(d_k)\cdot\gP\left(\sum_{m}E_{km}\otimes l_{km}\right)=\\
=\frac{1}{n^3}\left(\sum_{j}e^{j-k}\sigma^{-k+1}(l_{kj}^*\lambda_{j-1}^{-1}\lambda_{k-1}^{})\right)\cdot\sigma^{-k+1}(d_k)\cdot\left(\sum_{m}e^{k-m}\sigma^{-m+1}(l_{km})\right)=\\
=\frac{1}{n^3}\sum_{j,m}e^{j-k}\sigma^{-k+1}(l_{kj}^*\lambda_{j-1}^{-1}\lambda_{k-1}^{}d_k)e^{k-m}\sigma^{-m+1}(l_{km})=\\
=\frac{1}{n^3}\sum_{j,m}e^{j-k}e^{k-m}\sigma^{-m+1}(l_{kj}^*\lambda_{j-1}^{-1}\lambda_{k-1}^{}d_k)\sigma^{-m+1}(l_{km})=\frac{1}{n^3}\sum_{j,m}e^{j-m}\sigma^{-m+1}(l_{kj}^*\lambda_{j-1}^{-1}\lambda_{k-1}^{}l_{km}^{}d_k),
\end{gather*}
which proves (\ref{eq_main_cyc}). By Lemma \ref{lemma_PB_sigma_inv} we obtain $\sigma^{1-k}(d_k)\in P_B(\bL)$ and the assertion follows from (\ref{eq_sos_inv_r1}) and (\ref{eq_main_cyc}).
\end{proof}

\begin{remk}
From the proof we conlude that $\gP$ maps $\sum(M_n(\bL),\tau)^2$ onto $\sum\gA^2,$ that is, $\gP$ is even a strong conditional expectation.
\end{remk}

\begin{prop}\label{prop_P_B_nc_sos}
Every element $X\in P_B(\bL)\subseteq\gA$ is a non-commutative sum of squares in $\gA.$ In the case $n=2$ we have $P_B(\bL)\subseteq\sum\gA^2.$
\end{prop}
\begin{proof}
The assertion follows from the equalities $\lambda_k=e^{*k}e^k,\ k=1,\dots,n-1,$ and the fact that $P_B(\bL)$ is generated by $\lambda_1,\dots,\lambda_{n-1}.$
\end{proof}

\begin{prop}\label{prop_XY_pos_gA}
If $X,Y\in\gA_+$ and $XY=YX$, then $XY\in\gA_+.$ In particular, non-commutative sums of squares in $\gA$ belong to $\gA_+.$
\end{prop}
\begin{proof}
Combine Propositions \ref{prop_prod_com_pos} and \ref{prop_eps_*hom}.
\end{proof}

\begin{prop}\label{prop_pst_gA}
An element $x=x^*\in\gA$ is in $\gA_+$ if and only if it belongs to $QM_\gA(P_B(\bL)).$
\end{prop}
\begin{proof}
If $x\in\gA_+$, then we have $x\in(M_n(\bL),\tau)_+$ by Proposition \ref{prop_eps_*hom}. From Proposition \ref{prop_pst_alg_inv} we obtain $x\in QM_{M_n(\bL),\tau}(P_B^n(\bL))$ and from Proposition \ref{prop_gP_cyc_alg_main} we get $x\in QM_\gA(P_B(\bL)).$ On the other hand, since $P_B(\bL)\subseteq\gA_+$ we obtain $\gA_+\subseteq QM_\gA(P_B(\bL)).$
\end{proof}

Summarizing the preceding we obtain the following Posititvstellensatz. 
\begin{thm}\label{thm_pos_cyc_alg}
Let $x=x^*\in\gA$. Then $x\in \gA_+$ if and only $x$ is a non-commutative sum of squares. If $n=2$, then $\gA_+=\sum\gA^2.$
\end{thm}
\begin{proof}
Apply Propositions \ref{prop_P_B_nc_sos}, \ref{prop_pst_gA} and \ref{prop_XY_pos_gA}.
\end{proof}

\noindent{\bf Remarks.} 1. Let $\bL=\bR$. An ordering $p$ of $\bL$ is a $*$-ordering if $p$ contains $\gp(x^*x)$ for all $x\in\gA,$ where $\gp:\gA\to\bL$ is the canonical projection. As in the previous section one can say that {\it $\gA$ is "formally real" if and only if there exist a $*$-ordering on $\bL$}. For cyclic algebras this seems to be the most natural analogue of the notion of formally real field. In particular, it implies that 
\begin{gather}\label{eq_form_real}
\sum_ix_i^*x_i^{}=0\Rightarrow x_i=0.
\end{gather}

\noindent 2. Let us consider the complex $*$-algebra
$
\cA=\dC\langle a,b|\ a=a^*, b=b^*, ab=e^{\img\frac{2\pi}3}ba\rangle
$ and let $\gA$ denote the localization of $\cA$ by its center $Z(\cA)\simeq \dC[a^3,b^3].$ 
From the representation theory of this algebra studied in \cite{s5} it follows that $\cA$ has a faithful $*$-representation, say $\pi$, by unbounded $*$-reprentations. Therefore, (\ref{eq_form_real}) is valid for $\cA$ and hence for $\gA$. The algebra $\gA$ is a cyclic algebra associated with the extension $\dC(a^3,b)/\dC(a^3,b^3),$ but the involution of $\gA$ does not satisfy $\gA_k^*=\gA_{-k}.$ If we define $\cA_+$ by $\cA_+=\{x\in \cA:\pi(x)\geq 0\}$, it is natural to ask whether or not $x,y\in\cA_+$ and $xy=yx$ imply that $xy\in\cA_+.$

\noindent 3. All notions and results of this section remain valid in the following more general context. Suppose that $\bL/\bK$ is a Galois extension with group $G$ and $\Phi:G\times G\to\bL^\circ$ is a 2-cocycle. Let $\gA$ be the crossed product algebra $(\bL/\bK,\Phi),$ that is, $\gA=\oplus_{\sigma\in G}e_\sigma\cdot\bL$ as a right $\bL$-linear vector space with multiplication defined by
$$
(\sum_{\sigma\in G}e_\sigma c_\sigma)(\sum_{\rho\in G}e_\rho d_\rho):=\sum_{\sigma,\rho}e_{\sigma\rho}\Phi(\sigma,\rho)\rho(c_\sigma)d_\rho,
$$
see \cite{pi} for more details. Then $\gA$ is $G$-graded, where $\gA_\sigma=e_\sigma\cdot\bL.$ 

The results in \cite{cim} (see e.g. Theorem 10 therein) show that in order to conclude that (\ref{eq_form_real}) is valid it is natural to require $\gA_{\sigma}^{*}=\gA_{\sigma^{-1}}^{},\ \sigma\in G.$ 
%The notions of a $*$-ordering, positivity and all results of this section remain true.

\section{Examples.}\label{sect_exam}
Now we illustrate the results of previous two sections by  examples. All three examples in this section have the following in common: $\cA$ is a finitely generated complex $*$-algebra  and $\cA_+$ is defined by
$$
\cA_+:=\set{x\in\cA_h:\ \pi(x)\geq 0\ \mbox{for all finite-dimensional} *-\mbox{representations}\ \pi }.
$$ 
The center $Z(\cA)$ contains no zero divizors and the localization $\gA:=\cA(Z(\cA))^{-1}$ of $\cA$ by $Z(\cA)$  is a cyclic algebra of order $n$. Thus $\gA$ is a $*$-algebra over $\bK:=\Quot(Z(\cA)).$ In all three examples we  have the equality $\cA_+=\gA_+\cap\cA$ which allows us to apply the results from the previous sections.

Let $\Rep_n$ be the family of all $n$-dimensional irreducible $*$-representations of $\cA$ and set $n=2$ in the first example and $n=3$ in the second and third example.
Then the set $\Rep_n$ separates the elements of $\cA$ and an element $x\in \cA_h$ is in $\cA_+$ if (and only if) $\pi(x)\geq 0$ for all $\pi \in\Rep_n$.
 %finite-dimensional $*$-representations holds
%\begin{gather}\label{
%Further, $\cA$ has a faithful family $\Rep_n(\cA)$ of $n$-dimensional irreducible %$*$-representations and for the set $\cA_+$ of elements which are positive in all %finite-dimensional $*$-representations holds
%\begin{gather}\label{eq_pos_rep_n}
%\cA_+=\set{x\in\cA_h\ |\ \pi(x)\geq 0,\ \forall\pi\in\Rep_n(\cA)}.
%\end{gather}

\subsection{One-dimensional WCAR-algebra} In this subsection we set $$\cA:=\dC\langle x,x^*|x^*x+xx^*=1\rangle.$$ 
All irreducible $*$-representations of the $*$-algebra $\cA$ are of  dimensions $1$ or $2,$ see e.g. \cite{sav}. Every irreducible $2$-dimensional $*$-representation is of the form 
\begin{gather}\label{eq_rep_wcar}
x\mapsto\left(\begin{array}{ll} 0 & v \\ u & 0 \end{array}\right),\ x^*\mapsto\left(\begin{array}{ll} 0 & u \\ \overline{v} & 0 \end{array}\right),\ u\in\dR,\ v\in\dC,\ u^2+v\overline{v}=1.
\end{gather}
In the case $u^2=v\overline{v}$ each $*$-representation (\ref{eq_rep_wcar}) is  a sum of one-dimensional representations and all one-dimensional representations of $\cA$ are obtained in this way \cite{sav}. This implies that  $y\in\cA_h$ is in $\cA_+$ if and only if $y$ is positive in all $*$-representations of the form (\ref{eq_rep_wcar}). We  consider $u$ and $v$ as commuting variables satisfying $u^2+\overline{v}v=1,\ u=\overline{u}.$ Since the $*$-representations (\ref{eq_rep_wcar})  separate the elements of $\cA$, they define a $*$-embedding $\cA\hookrightarrow M_2(\dC[u,v]).$ Via this embedding $\cA$ is identified with the $*$-algebra of matrices of the  form
\begin{gather}\label{eq_embed_wcar}
\left(\begin{array}{ll} P_1+u^2P_2 & vP_3+uP_4 \\ uP_3+\overline{v}P_4 & P_1+v\overline{v}P_2 \end{array}\right),\ P_i\in\dC[uv,u\overline{v}].
\end{gather}
%and (\ref{eq_pos_rep_n}) holds with $n=2.$ 
The center $Z(\cA)$ is isomorphic to $\dC[uv,u\overline{v}]$ and it is an integral domain. We denote by $\gA$ the algebra $\cA(Z(\cA))^{-1},$ by $\bK$ the field $\Quot(Z(\cA))=\dC(uv,u\overline{v}),$ and by $\bL$ the field $\dC(u^2,uv,u\overline{v}).$% We identify $\gA$ with the subalgebra of matrices in $M_2(\dC(u,v,\overline{v}))=M_2(\bL_1).$

The Galois extension $\bL/\bK$ is defined by the $\bK$-automorphism $$\sigma:u^2\mapsto v\overline{v}=1-u^2,$$ of $\bL$ which is of order 2. Let $e$ be the generator $x.$ Then $e^2=\left(\begin{array}{ll} uv & 0 \\ 0 & uv \end{array}\right)=uv\cdot \Un{\gA}.$ Then one easily checks  that $\gA$ is the cyclic algebra $(\bL/\bK,\sigma,uv)$ which satisfies assumptions of Section \ref{sect_sos_cyc}. Let $\gA_+$ be as in Section \ref{sect_sos_cyc}. 

%Let $\cA_+$ denote the set of $x=x^*\in\cA$ which are positive in all $*$-representations 
\begin{lemma}\label{prop_wcar_A+_gA+}
$\cA_+=\gA_+\cap\cA.$
\end{lemma}
\begin{proof}
We only  sketch  the proof. If we take in (\ref{eq_embed_wcar}) $P_i\in\bK=\dC(uv,u\overline{v})$ we obtain a $*$-embedding $\gA\to M_2(\dC(u,v,\overline{v})).$ We define a mapping $\phi:\gA\to M_2(\bL)$ by
\begin{gather}\label{eq_embed_wcar_gA}
\left(\begin{array}{ll} P_1+u^2P_2 & vP_3+uP_4 \\ uP_3+\overline{v}P_4 & P_1+v\overline{v}P_2 \end{array}\right)\mapsto\left(\begin{array}{ll} P_1+u^2P_2 & uvP_3+u^2P_4 \\ P_3+\frac{u\overline{v}}{u^2}P_4 & P_1+v\overline{v}P_2 \end{array}\right),\ P_i\in\bL.
\end{gather}
Let $B=\diag{1,u^2}\in M_2(\bL)$ and let $\tau$ be the corresponding involution on $M_2(\bL),$ that is $\tau(X)=B^{-1}X^*B,\ X\in M_2(\bL).$ A direct computation shows that $\phi$ is a $*$-embedding of $\gA$ into $(M_2(\bL),\tau).$ As in Section \ref{sect_sos_matr_field}, let $P_B(\bL)$ be the preordering in $\bL$ generated by $u^2.$

Take $y=y^*\in\cA_+.$ Let $\Delta_1,\Delta_2,\Delta_3\in\dC[u,v,\overline{v}]$ be the principal minors of $y.$ Then $\Delta_i\in\dC[u,v,\overline{v}]_+\subseteq\dC(u,v,\overline{v})_+.$ From  (\ref{eq_embed_wcar}) it follows that $\Delta_i\in\bL.$ From the Krivine-Stengle theorem (see \cite{mar}) we conclude that $\dC(u,v,\overline{v})_+=\sum\dC(u,v,\overline{v})^2.$ A simple computation shows that $\sum\dC(u,v,\overline{v})^2\cap\bL=\sum\bL^2+u^2\sum\bL^2=P_B(\bL).$ Hence $\Delta_i\in P_B(\bL).$ It follows from (\ref{eq_embed_wcar_gA}) that the corresponding leading minors of $\phi(y)$ are also equal to $\Delta_i.$ Since $\Delta_i\in P_B(\bL)$ by Lemma \ref{lemma_pos_in_Mn_K}, we get $\phi(y)\in (M_2(\bL),\tau)_+$ and Proposition \ref{prop_eps_*hom} implies that $y\in\gA_+.$ 

A similar reasoning shows that $y\in\gA_+\cap\cA$ implies $y\in\cA_+.$
\end{proof}

% \begin{proof}
% We only give a sketch of a proof which uses theory of induced $*$-representations developed in \cite{ss}. Algebra $\cA$ has a natural $\cyc{2}$-grading $\cA_0\oplus\cA_1$ defined by $\cA_i=\gA_i\cap\cA.$ Let $\cB=\cA_0,$ and $\gp:\gA\to\gA_0$ be the canonical projection, thus defining a conditional expectation from $\cA$ onto $\cB$ denoted by $\gp$ as well. As in \cite{ss}, let $\cBp$ be the set of characters $\chi:\cB\to\dC$ such that $\chi(\sum\cA^2\cap\cB)\geq 0.$ Since all irreducible representations of $\cA$ are induced, an element $x=x^*\in\cA$ belongs to $\cA_+$ if and only if $\pi_\chi(x)\geq 0$ for all $\chi\in\cBp,$ where $\pi_\chi=\Ind(\chi).$ By construction of $\Ind(\chi)$ it is equivalent to $\chi(\gp(y^*xy))\geq 0$ for all $y\in\cA$ and all $\chi\in\cBp.$ A character $\chi\in\cBd$ belongs to $\cBp$ if and only if $\chi(a^*a)\geq 0.$ Thus, $\chi(y^*xy)\geq 0$ for all $\chi\in\cBp$ if and only if $y^*xy$ is positive in all $*$-orderings of $\bL.$ By Definition \ref{defn_pos_in_gA} it is equivalent to $x\in\gA_+.$
% \end{proof}

Combining Lemma \ref{prop_wcar_A+_gA+} with Theorem \ref{thm_pos_cyc_alg} we obtain the following

\begin{thm}
An element $y=y^*\in\cA$ is in $\cA_+$ (that is, $y$ is positive in all finite-dimensional $*$-representations of $\cA$) if and only if $y\in\sum\gA^2$, or equivalently, there exists a $c\in Z(\cA),\ c\neq 0$, such that $c^*c\cdot y\in\sum\cA^2.$
\end{thm}

\begin{remk}
The preceding theorem can be also obtained by using Theorem 5.4 and Corollary 5.5 in \cite{ps}, since $\gA$ is a quaternion algebra. 
\end{remk}

\subsection{An algebra related to $\widetilde{E_6}$.} In this subsection $\cA$ is a $*$-algebra related to the extended Dynkin diagram $\widetilde{E_6}$ (see \cite{m} and the references therein), that is,  
\begin{gather*}
\cA:=\mathbb{C}\langle a_1,a_2,a_3|a_1+a_2+a_3=0, a_i^*=a_i^{}=a_i^3\rangle.
\end{gather*}

We can generate $\cA$ by the so-called \textit{centered} element $x$ (see \cite{m}) which  is defined by relations $a_i=\eps^ix+\eps^{-i}x,\ i=1,2,3,$ where $\varepsilon=e^{\img\frac{2\pi}{3}}.$ Then we have the following lemma \cite{m}:
\begin{lemma}
If $x$ and $x^*$ are taken as generators, then $\cA$ has the form 
\begin{gather*}
\cA=\dC\langle x,x^*|x^3+x^{*3}=0,\ x^2x^*+xx^*x+x^*x^2-x=0,\\ x^{*2}x+x^*xx^*+xx^{*2}-x^*=0\rangle
\end{gather*}
\end{lemma}

The $*$-algebra $\cA$ has the following family of $3$-dimensional $*$-representations 
\begin{gather}\label{eq_rep_E6}
x\mapsto\left(\begin{array}{lll} 0 & 0 & \img v_3\\ v_1 & 0 & 0\\ 0 & v_2 & 0 \end{array}\right),\ 
x^*\mapsto\left(\begin{array}{lll}\ \ \ 0 & v_1 & 0\\ \ \ \ 0 & 0 & v_2\\ -\img v_3 & 0 & 0 \end{array}\right)
\end{gather}
where $v_1,v_2,v_3\in\dR$ satisfying $v_1^2+v_2^2+v_3^2=1.$ It can be shown that the representations (\ref{eq_rep_E6}) form a separating family, so they define a $*$-embedding $\cA\hookrightarrow M_3(\dC[v_1,v_2,v_3]).$ In this manner $\cA$ is identified with the $*$-algebra of matrices of the following form:
\begin{gather*}
\left(\begin{array}{lll} \ \ \ \ P & v_1\cdot Q & \img v_3\cdot R\\ \ \ v_1\cdot \sigma(R) & \sigma(P) & v_2\cdot \sigma(Q)\\ -\img v_3\cdot \sigma^2(Q) & v_2\cdot \sigma^2(R) & \sigma^2(P)\end{array}\right),\ \mbox{where}\ P,Q,R\in\dC[v_1^2,v_2^2,v_3^2,v_1v_2v_3].
\end{gather*}
Let $\sigma$ be the automorphism of $\cB:=\dC[v_1^2,v_2^2,v_3^2,v_1v_2v_3]$ defined by
$$
v_1^2\mapsto v_2^2,\ v_2^2\mapsto v_3^2,\ v_3^2\mapsto v_1^2,\ v_1v_2v_3\mapsto v_1v_2v_3.
$$

Then the center $Z(\cA)$ is isomorphic to the stable subalgebra of $\cB$  under the automorphism $\sigma$ and it is an integral domain. Let $\gA$ denote the algebra $\cA(Z(\cA))^{-1},$ $\bK$ the field $\Quot(Z(\cA)),$ and $\bL$ the field $\dC(v_1^2,v_2^2,v_3^2,v_1v_2v_3).$ 

The automorphism $\sigma$ extends to $\bL$ and defines the Galois extension $\bL/\bK,$ which is a cyclic extension, since $\sigma^3=\id.$ Let $e$ be  the generator $x.$ Then 
$e^3=\img v_1v_2v_3\cdot \Un{\gA}.$ One checks by some direct computations that $\gA=\bL+\bL e+\bL e^2$ is the cyclic algebra $(\bL/\bK,\sigma,\img v_1v_2v_3).$ 

As in the previous subsection the following lemma can be proved.
\begin{lemma}\label{prop17}
$\cA_+=\gA_+\cap\cA.$
\end{lemma}

Combining Lemma \ref{prop17} and Theorem \ref{thm_pos_cyc_alg} it follows that an element $y=y^*\in\cA$ is in $\cA_+$ if and only if it is in $\sum_{nc}\cA^2$, that is, $y$ is a non-commutative sum of squares in $\gA.$ In fact, the following stronger result  holds.
\begin{thm}
Let $y=y^*\in\cA$. Then $y\in\cA_+$ if and only if $y\in\sum\gA^2$, or equivalently, there exists $c\in Z(\cA),\ c\neq 0,$ such that $c^*c\cdot y\in\sum\cA^2.$
\end{thm}
\begin{proof}
By Proposition \ref{prop_pst_gA} we have $\gA_+=QM_\gA(P_B(\bL)).$ Hence $\gA_+$ is the quadratic module generated by $e^*e\cdot e^{*2}e^2$ or equivalently by $e^*e\cdot ee^*.$ Since $e=x$, we have
$e=e^*e^2+e^2e^*+ee^*e.$ Multiplying this equation by $e^*$ from the left and by $ee^*$ from the right and remembering that $e^3\in Z(\gA),\ ee^{*2}e=e^*e^2e^*,$ we derive
\begin{gather*}
e^*e\cdot ee^*=e^{*2}e^3e^*+e^*e^2e^*ee^*+e^*ee^*e^2e^*=\\
=e^{*3}e^3+ee^{*2}e^2e^*+e^*e^2e^{*2}e=e^{*3}e^3+(e^2e^*)^*e^2e^*+(e^{*2}e)^*e^{*2}e\in\sum\gA^2.
\end{gather*}
Thus $\gA_+=QM(e^*e\cdot ee^*)=\sum\gA^2.$
\end{proof}

\subsection{A counterexample to a question of Procesi and Schacher}\label{subsect_ctrexam} Procesi and Schacher  \cite{ps} asked if the denominator-free Positivstellensatz holds in central simple algebras (CSA) with involution. Recently Klep and Unger \cite{klep} gave a nice counterexample. We now provide another counterexample which is a cyclic algebra. 

Let $\dR[x,y,z]$ be the $*$-algebra of polynomials in three real variables and $\cA$ be the $*$-subalgebra of $M_3(\dR[x,y,z])$ generated by the identity and the matrices
\begin{gather}
X=\left(\begin{array}{lll} 0 & 0 & z\\ x & 0 & 0\\ 0 & y & 0 \end{array}\right),\ 
X^*=\left(\begin{array}{lll}0 & x & 0\\ 0 & 0 & y\\ z & 0 & 0 \end{array}\right).
\end{gather}

As in the previous subsection, let $\gA\subseteq M_3(\dR(x,y,z))$ denote the localization of $\cA$ by its center. We consider $\gA$ with the natural involution inherited from $\cA.$ 

Let $\bL=\dR(x^2,y^2,z^2,xyz)$ and let $\sigma$ be the automorphism of $\bL$ of order three defined by $$\sigma(x^2)=y^2,\ \sigma(y^2)=z^2,\ \sigma(z^2)=x^2,\ \sigma(xyz)=xyz.$$ We denote by $\bK\subset\bL$ the stable subfield  under $\sigma.$ We identify $\bL$ with a subfield of $M_3(\dR(x,y,z))$ via the embedding
\begin{gather}\label{eq_aux_1}
\bL\ni l\mapsto\diag{l,\sigma(l),\sigma^2(l)}. 
\end{gather}

By direct computations one checks that 
\begin{gather}\label{eq_ctrex_sum}
\gA=\bL\oplus\bL X\oplus\bL X^*=\bL\oplus\bL X\oplus\bL X^2.
\end{gather}

\begin{lemma}
$\gA$ is isomorphic to the cyclic algebra $(\bL/\bK,\sigma,xyz).$
\end{lemma}
\begin{proof}
This follows from (\ref{eq_ctrex_sum}) and the following equalities 
\begin{gather*}
\diag{l,\sigma(l),\sigma^2(l)}\cdot X=X\cdot\diag{\sigma(l),\sigma^2(l),l}\in M_3(\dR(x,y,z)),\\
X^*=\diag{\frac{x^2}{xyz},\frac{y^2}{xyz},\frac{z^2}{xyz}}\cdot X^2\ \mbox{and}\ X^3=xyz\cdot I.
\end{gather*}
\end{proof}

\begin{lemma}\label{lemma_y_11}
Let $Y=(y_{ij})\in\gA$ be an element of $\sum\gA^2.$ Then $y_{11}\in\sum\bL^2+x^2\cdot\sum\bL^2+z^2\cdot\sum\bL^2.$
\end{lemma}
\begin{proof}
Take $Z=l_0+l_1X+l_2X^*\in\gA$ and let $Y=ZZ^*.$ Then $y_{11}=l_0^{}l_0^*+z^2l_1^{}l_1^*+x^2l_2^{}l_2^*.$
\end{proof}

The following lemma occurs also  in \cite{klep}. We include an elementary proof.
\begin{lemma}\label{lemma_aux_1}
Let 
\begin{gather}\label{eq_aux_2}
s_0xy=s_1+s_2x+s_3y,\ s_i\in\sum\dR(x,y)^2,\ i=0,1,2,3.
\end{gather}
Then $s_i=0$ for $ i=0,1,2,3.$
\end{lemma}
\begin{proof}
After multiplying (\ref{eq_aux_2}) by a common denominator we are reduced to the case $s_i\in\sum\dR[x,y]^2.$ Dividing both sides by a power of $x^2$ we can assume that $x^2$ does not divide all summands in $s_i.$

Setting $x=0,\ y>0,$ we get $0=s_1(0,y)+ys_3(0,y)$ which implies $s_1(0,y)=s_3(0,y)=0$, so that $s_1$ and $s_3$ are divisible by $x.$ Thus, each summand in $s_1,s_3$ is divisible by $x,$ hence by $x^2.$ That is, $s_1=s_1'x^2,\ s_3=s_3'x^2,$ where $s_1',s_3'\in\sum\dR[x,y]^2.$ Cancelling $x$ in (\ref{eq_aux_2}) we get
\begin{gather}\label{eq_aux_3}
s_0y=s_1'x+s_2+s_3'xy.
\end{gather}
Setting $x=0,\ y<0,$ in (\ref{eq_aux_3}) we obtain  in similar manner $s_0=x^2s_0',\ s_2=x^2s_2'.$ Thus all summands in $s_i,\ i=0,1,2,3,$ are divisible by $x^2$ which is a contradiction.
\end{proof}

\begin{lemma}\label{lemma_ctrex_impos}
Suppose that we have an equality 
\begin{gather}\label{eq_aux_4}
s_0xy+s_1z=s_2+s_3x+s_4y+s_5xz+s_6yz+s_7xyz,\ s_i\in\sum\dR(x,y,z)^2,\ i=0,\dots,7.
\end{gather}
Then  $s_i=0$ for all $i=0,\dots,7.$
\end{lemma}
\begin{proof}
The proof uses the same reasoning as in the proof of Lemma \ref{lemma_aux_1}. We consider only the case $s_i\in\sum\dR[x,y,z]^2$ and assume that $z^2$ does not divide all summands in $s_i.$ 

Setting $z=0$ in (\ref{eq_aux_4}) we get 
\begin{gather}
s_0(x,y,0)xy=s_2(x,y,0)+s_3(x,y,0)x+s_4(x,y,0)y.
\end{gather}
Lemma \ref{lemma_aux_1} implies that $s_0,s_2,s_3,s_4$ are divisibe by $z^2,$ i.e. $s_0=s_0'z^2,s_2=s_2'z^2,s_3=s_3'z^2,s_4=s_4'z^2$, where $ s_i'\in\sum\dR[x,y,z]^2,\ i=0,2,3,4.$ Dividing both sides of (\ref{eq_aux_4}) by $z$ we obtain
\begin{gather}\label{eq_aux_5}
s_0'xyz+s_1=s_2'z+s_3'zx+s_4'zy+s_5x+s_6y+s_7xy,\ s_i,s_j'\in\sum\dR(x,y,z)^2.
\end{gather}
Seeting $z=0$ and dividing both sides by $xy$ in (\ref{eq_aux_5}) we derive 
$$
\frac{s_1(x,y,0)}{(xy)^2}xy=\frac{s_5(x,y,0)}{y^2}s_5y+\frac{s_6(x,y,0)}{x^2}s_6x+s_7(x,y,0).
$$
Lemma \ref{lemma_aux_1} implies that $s_1(x,y,0)=s_5(x,y,0)=s_6(x,y,0)=s_7(x,y,0)=0.$ Hence all elements $s_i$ are divisible by $z^2$ which is a contradiction. 
\end{proof}

Then an element $(y_{ij})\in\cA\subset M_3(\dR[x,y,z])$ is in $\cA_+$ if and only if the matrix $(y_{ij}(x,y,z))$ is positive semi-definite for all $(x,y,z)\in\dR^3.$ As in the previous two examples, one can prove that $\cA_+=\cA\cap\gA_+,$ where $\gA_+$ is defined as in Section \ref{sect_sos_cyc}.

\begin{prop}
The matrix $Y=X^*X^2X^*\in\cA$ is a positive element of $\cA$ such that $Y\notin\sum\gA^2,$ or equivalently, there is no element $c\in Z(\cA),\ c\neq 0$ such that $c^*c\cdot Y\in\sum\cA^2.$
\end{prop}
\begin{proof}
The element $Y$ is equal to $\diag{x^2z^2,y^2x^2,z^2y^2},$ hence $Y\in\cA_+.$ Assume to the contrary that $Y\in\sum\gA^2.$ Then by Lemma \ref{lemma_y_11}  there is an equality 
\begin{gather}\label{eq_ctrex_impos}
s_0x^2z^2=s_1+s_2x^2+s_3z^2,\ \mbox{where}\ \ s_i\in\sum\in\dR(x^2,y^2,z^2,xyz)^2,\  s_0=1.
\end{gather}

Each element $t\in\dR(x^2,y^2,z^2,xyz)$ can be written as $t_1+xyz t_2,\ t_i\in\dR(x^2,y^2,z^2).$ Hence $s_i=r_i+q_ix^2y^2z^2,\ r_i,q_i\in\sum\dR(x^2,y^2,z^2)^2.$ Applying this to (\ref{eq_ctrex_impos}) we obtain 
\begin{gather*}
r_0x^2z^2+(q_0(x^2z^2)^2)y^2=r_1+q_1x^2y^2z^2+r_2x^2+(x^4q_2)y^2z^2+r_3z^2+(z^4q_3)x^2y^2.
\end{gather*}
Applying Lemma \ref{lemma_ctrex_impos} to $\dR(x^2,y^2,z^2)$ we get $r_i=q_i=0,$ hence $s_i=0$ for all $i$. This is a contradiction, since $s_0\neq 0$.

\end{proof}

\section{Some problems}\label{openproblems}
% 
% Let $G$ be a discrete group and $H\subseteq G$ be a subgroup. Assume in addition that index $[G:H]=n$ is finite. Fix one element $g_i,\ i=1,\dots,n$ in each left coset of $G$ by $H.$ So that $G=\cup_{i=1}^n g_iH.$ Let $\cB=\dC[H],\ \cA=\dC[G]$ and $e_i=g_i.$ Then conditions (H1)-(H5) in the Appendix are satisfied and there exists a strong conditional expectation from $M_n(\dC[H])$ onto $\dC[G].$
% 
% {\bf Problem.} {\it Does Positivstellensatz I holds for the algebra $M_n(\dC[F_m]),$ where $F_m$ is the free group of $m$ generators.}
% 
% A positive answer to Problem 1 would imply the following proposition: {\it If $G$ be a group having a free subgroup of a finite index, then $psd=sos$ holds in $\dC[G].$}
% 
% {\bf Problem.} {\it Let $(M_n(\dL),\tau_B)$ be the algebra of matrices over a real field $\dL$ and invlution $\tau_B$ defined by a non-sigular matrix $B.$ Assume that $\tau_B$ is real. Under what conditions on $B$ we have $psd=sos$ in $(M_n(\dL),\tau_B)$}
% 
% {\bf Problem.} {\it Find an example of a $*$-algebra $\cA$, such that $psd=sos$ holds in $\cA$ and does not hold in $\cA.$}

Let $\cW(d)$ be the Weyl algebra, that is, $\cW(d)$ is the unital $*$-algebra with self-adjoint generators $p_1,\dots,p_d,q_1,\dots,q_d$ and defining relations
\begin{gather*}
p_jp_k=p_kp_j,\ q_jq_k=q_kq_j,\ p_jq_k=q_kp_j\ \mbox{for}\ j,k=1\cdots,d, j\neq k,\\ 
p_kq_k-q_kp_k=-i\ \mbox{for}\ k=1,\dots, d.
\end{gather*}
There is a distingushed faithful $*$-representation $\pi_0$ of $\cW(d)$ on the Schwartz space $\cS(\dR^d)$, called the {\it Schr\"odinger representation}, defined by 
$$
(\pi_0(p_k)f)(t)=\frac{\partial}{\partial t_k}f(t),~(\pi_0(q_k)f)(t)= t_kf(t),~~k=1,\cdots,d,~ f\in \cS(\dR^d).$$ 
Define $$\cW(d)_+:= \{x \in \cW(d): \langle \pi_0(x)f,f \rangle \geq 0 ~{for}~ f\in \cS(\dR^d)\}.$$

\smallskip
\noindent{\bf Problem 1:} {\it Does a Positivstellensatz of type II hold for the algebra $\cW(\dR^d)$, that is, given $x\in \cW(d)_+$, does there exist $c\in \cW(d)$, $c\neq 0$, such that $cxc^*\in \sum \cW(d)^2$? If yes, can $c$ be chosen such that the kernel of the of operator $\ov{\pi_0(c)}$ is contained in the kernel of $\ov{\pi_0(x)}$?}

\smallskip
Strict Positivstellens\"atze and results supporting this question were proved in \cite{sweyl} and in \cite{sfrac}.

\bigskip
There is a similar problem for enveloping algebras of finite dimensional Lie algebras. Let $\gog$ be a real Lie algebra. Then the complex universal enveloping algebra $\cE(\gog)$ of $\gog$ is a complex unital $*$-algebra with involution determined by $x^*:=-x$ for $x\in\gog.$
Let $G$ denote the connected simply connected Lie group which $\gog$ as its Lie algebra and let $\widehat{G}$ be the set of unitary equivalence classes of irreducible unitary representations of $G.$ For each $U\in\widehat{G}$ there is an associated $*$-representation $dU$ of $\cE(\gog)$ with domain $\cD^\infty(U)$, see \cite{s1}, Chapter 10. Now let $$\cE(\gog)_+:= \set{x \in \cE(\gog): \langle dU(x)f,f \rangle \geq 0\ \mbox{for}\ f\in\cD^\infty(U),\ U\in\widehat{G}}.$$ The counterpart of Problem 1 for enveloping algebras is the following

\smallskip
\noindent{\bf Problem 2:} {\it Is a Positivstellensatz of type II true for $\cE(\gog)$, that is, given $x\in \cE(\gog)_+$, does there exist an element $c\in \cE(\gog),\ c\neq 0,$ such that $cxc^*\in \sum \cE(\gog)^2$?}

\bigskip
In both types III and IV  commuting positive elements occur. While the product of two commuting positive bounded operators on a Hilbert space is always positive (Lemma \ref{lemma_prod_com}), there are examples of  commuting positive symmetric operators on a unitary space for which the product is no longer positive. 
It seems to be unknown whether or not the latter can happen in the Schr\"odinger representation of the Weyl algebra. 

\noindent{\bf Problem 3:} {\it Suppose that $a,b\in\cW(d)_+$ and $ab=ba.$ Is it true that $ab\in\cW(d)_+$?}

\smallskip
An affirmative answer would imply that all elements of the minimal non-commutative preordering in the Weyl algebra are indeed positive elements.

\bigskip
Let $\cA$ be a unital $*$-algebra and $\cR$ a separating family of $*$-representations of $\cA.$ If $\pi$ is a $*$-representation of $\cA$ on $\cD$, there is a unique $*$-representation $\pi_n$ of the matrix $*$-algebra $M_n(\cA)$ on $\cD_n=\cD\oplus\cdots\oplus\cD$ ($n$ times) defined by $\pi_n((a_{kl}):=(\pi(a_{kl})).$

Let $\cA_+:=\{a\in \cA:\pi(a)\geq 0\ \mbox{for}\ \pi\in\cR\}$ and $M_n(\cA)_+:=\{A\in \cA:\pi_n(A)\geq 0~~{\rm for}~~\pi \in \cR\}$.

\smallskip
\noindent{\bf Problem 4:} {\it Suppose that a Positivstellensatz of type I hold for $\cA.$ Does it hold also for the matrix algebra $M_n(\cA)$?}

In particular, Problem 4 is open and important when $\cA$ is the commutative real $*$-algebra of all polynomials on the $2$-sphere $S^2=\{(x,y,z)\in \dR^3: x^2+y^2+z^2=1\}$ and $\cR$ is the set of all point evaluations $\pi_t(p)=p(t),\ t\in S^2$. It follows from results in \cite{sd} that each nonnegative polynomial on $S^2$ is a sum of squares of polynomials, that is, a Positivstellensatz of type I is valid for $\cA$. By the results of Section \ref{sect_nc_matr} there is a Positivstellensatz of type II for the matrix algebra $M_n(\cA)$. The question is whether or not a {\it denominator free} Positivstellensatz holds for $M_n(\cA)$.

An affirmative answer to Problem $4$ in the latter case would yield a number of other interesting results. First, using a similar conditional expectation as in Section \ref{sect_sos_interval} it would follow that \textit{a self-adjoint matrix polynomial $F(x,y)\in M_n(\dC[x,y])$ is positive semi-definite on the unit disc $\set{x^2+y^2\leq 1}$ if and only if $F\in\sum M_n(\dC[x,y])^2+(1-x^2-y^2)\sum M_n(\dC[x,y])^2.$} 

Secondly, it would imply that a Positivstellensatz of type I holds for a number of algebras which can be embedded into $M_n(\dC[S^2]).$ For example, consider the following "non-sommutative sphere"
$$
\cA=\dC\langle x_1,x_2,x_3|x_1^2+x_2^2+x_3^2=I,\ x_ix_j=-x_jx_i,\ x_i^*=x_i,\ i=1,2,3\rangle.
$$
It follows from the description of irreducible representations of $\cA$ (see e.g. \cite{osam}, p.110) that $\cA$ can be embedded into $M_2(\dC[S^2])$ such that $\cA_+\subseteq M_2(\dC[S^2])_+.$ Further, there exists a strong conditional expectation $\gP:M_2(\dC[S^2])\to\cA.$ Hence an affirmative answer to Problem 4 for the algebra of polynomials on the sphere $S^2$ would give a Positivstellensatz of type I for $\cA.$

\bigskip
Let $G$ be a discrete group and let $\dC[G]_+$ denote the set of elements in the group algebra which are positive in all $*$-representations. 

\smallskip
\noindent{\bf Problem 5:} {\it For which discrete groups $G$ a Positivstellensatz of type I holds, that is, when is $\dC[G]_+=\sum\dC[G]^2$?}

First let $G=\dZ^n$. Then $\dC[\dZ^n]$ is isomorphic to the coordinate ring of the $n$-torus $\dT^n.$ or equivalently, to the $*$-algebra of  trigonometric polynomials in $n$ variables. Therefore the answer is affirmative for $n=1$ by the classical Fejer-Riesz theorem and for $n=2$ by the results in \cite{sd} and it is negative for $n\geq 3$ (\cite{ru}, see e.g. \cite{dr}). Using these facts it is easily shown that for an abelian group $G$ the equality  $\dC[G]_+=\sum\dC[G]^2$ holds if and only if the torsion-free component of $G$ is either $\dZ$ or $\dZ^2.$ 

Using similar techniques as in to \cite{hpm} one can  show $\dC[G]_+=\sum\dC[G]^2$ if $G$ is a free group. 

Let $H\subseteq G$ be a subgroup of a finite index $[G:H]=n.$ Then there is a natural construction of an embedding $\dC[G]\subseteq M_n(\dC[H])$ and of a strong conditional expectation $\gP:M_n(\dC[G])\to\dC[H]$ similar to the  one  in Section \ref{sect_sos_cross}. Recall that $G$ is called virtually free (resp. cyclic) if $G$ contains a free (resp. cyclic) group of a finite index. In  view of the Problem 4 we  conjecture that {\it $\dC[G]_+=\sum\dC[G]^2$ holds for virtually free and virtually cyclic groups of rank 2.} In the case when $\dZ\subseteq G$ and $[G:\dZ]<\infty$  the equality $\dC[G]_+=\sum\dC[G]^2$ can be obtained from Proposition \ref{prop_fej_rsz}.

\bigskip
\noindent{\bf Acknowledgements.} We would like to thank Andreas Thom and Tim Netzer for fruitful discussions, in particular on the subject of the last section.

\bibliographystyle{amsalpha}

\end{document}